  \documentclass[10pt]{amsart}
\usepackage{amsmath}
\usepackage{amsthm}
\usepackage{amscd}
\usepackage{amssymb}
\usepackage{amsfonts}
\usepackage{graphics}
\usepackage{color}
\usepackage[USenglish]{babel}
\usepackage{enumitem}

\date{}

\newlength{\defbaselineskip}
\long\def\salta#1{\relax}

\numberwithin{equation}{section}

\theoremstyle{plain}
\newtheorem{theorem}{Theorem}[section]
\newtheorem{proposition}[theorem]{Proposition}

\theoremstyle{definition}
\newtheorem{definition}[theorem]{Definition}

\newtheorem{remark}[theorem]{Remark}

\theoremstyle{remark}

\renewcommand{\theequation}{\thesection.\arabic{equation}}

\def\eps{\varepsilon}

\def\dys{\displaystyle}

\def\oeps{\Omega^\eps}

\def\t1p0{T^{1,p}_{0}(\Omega)}

\def\m2{M^{\frac{N(p-1)}{N-1}}(\Omega)}

\def\into{\int_{\Omega}}

\def\w-1p'{W^{-1,p'}(\Omega)}
\def\pw-1p'u{L^{p'}(0,1;W^{-1,p'}(\Omega))}

\def\dys{\displaystyle}

\def\lp'n{(L^{p'}(\Omega))^{N}}

\def\supp{\text{\rm{supp}}}

\def\oeps{\Omega^\eps}

\usepackage{scalerel}
\usepackage{stackengine}

\DeclareMathSymbol\mmho              {\mathord}{AMSb}{"5D}
\newcommand\reallywidetilde[1]{%
  \savestack{\tmpbox}{\stretchto{%
      \scaleto{%
        \scalerel*[\widthof{\ensuremath{#1}}]{\kern-1pt\raisebox{-.0005em}{$\mmho$}\kern-1pt}%
        {\rule[-\textheight/2]{1ex}{\textheight}}
      }{.8\textheight}
    }{2ex}}%
   \stackon[-1.02ex]{$#1$}{\tmpbox}%
}

\numberwithin{equation}{section}

\title[Homogenization of a strongly singular  problem]{ Homogenization\\ of a Dirichlet semilinear elliptic problem\\ with a strong singularity at $u=0$\\ in a domain with many small holes}

\author[D. Giachetti]{Daniela Giachetti}

\address{Daniela Giachetti \newline Dipartimento di Scienze di Base e Applicate per l'Ingegneria \newline Facolt\`a di ingegneria Civile e Industriale, Sapienza Universit\`a di Roma  \newline Via Scarpa 16, 00161 Roma, Italy}
\email{{\tt daniela.giachetti@sbai.uniroma1.it}}

\author[P.J. Mart\'inez-Aparicio]{Pedro J. Mart\'inez-Aparicio}

\address{Pedro J. Mart\'inez-Aparicio \newline Departamento de Matem{\'a}tica Aplicada y Estad\'istica \newline
Universidad Polit{\'e}cnica de Cartagena \newline Paseo Alfonso XIII 52, 30202 Cartagena (Murcia), Spain}
\email{{\tt pedroj.martinez@upct.es}}

\author[F. Murat]{Fran\c cois Murat}
\address{Fran\c cois Murat \newline Laboratoire Jacques-Louis Lions \newline Universit\'e Pierre et Marie Curie et C\MakeLowercase{nrs} \newline Bo\^ite Courrier 187, 75252 Paris Cedex 05, France}
\email{{\tt murat@ann.jussieu.fr}}

\keywords{homogenization, perforated domains, strange term, semilinear elliptic problem, strong singularity at $u = 0$
\\
\indent 2010 {\it Mathematics Subject Classification.} 35B25, 35B27, 35J25, 35J67}

\date{}

\begin{document}
\maketitle

{\bf \centerline {Version May 26, 2017}}
{\bf \centerline{Submitted for publication}}

\begin{abstract}
 In the present paper we perform the homogenization of the semilinear elliptic problem
  \begin{equation*}
\begin{cases}
u^\eps \geq 0 & \mbox{in} \; \oeps,\\
\displaystyle - div \,A(x) D u^\eps  = F(x,u^\eps) & \mbox{in} \; \oeps,\\
u^\eps = 0 & \mbox{on} \; \partial \oeps.\\
\end{cases} 
\end{equation*}
In this problem $F(x,s)$ is a Carath\'eodory function such that\break $0 \leq F(x,s) \leq h(x)/\Gamma(s)$ a.e. $x\in\Omega$ for every $s > 0$, with $h$ in some $L^r(\Omega)$ and $\Gamma$ a $C^1([0, +\infty[)$ function such that $\Gamma(0) = 0$ and $\Gamma'(s) > 0$ for every $s > 0$. On the other hand the open sets $\oeps$ are obtained by removing many small holes from a fixed open set $\Omega$ in such a way that a ``strange term" $\mu u^0$ appears in the limit equation in the case where the function $F(x,s)$ depends only on $x$.

    We already treated this problem in the case of a ``mild singularity", namely in the case where the function $F(x,s)$ satisfies  $0 \leq F(x,s) \leq h(x) (\frac 1s + 1)$. In this case the solution $u^\eps$ to the problem belongs to $H^1_0 (\oeps)$ and its definition is a ``natural" and rather usual one.
    
     In the general case where $F(x,s)$ exhibits a ``strong singularity" at $u = 0$, which is the purpose of the present paper, the solution $u^\eps$ to the problem only belongs to $H_{\mbox{\tiny{loc}}}^1(\oeps)$ but in general does not belongs to $H^1_0 (\oeps)$ any more, even if $u^\eps$ vanishes on $\partial\oeps$ in some sense. Therefore we introduced a new notion of solution (in the spirit of the solutions defined by transposition) for problems with a strong singularity. This definition allowed us to obtain existence, stability and uniqueness results.

   In the present paper, using this definition, we perform the homogenization of the above semilinear problem and we prove that in the homogenized problem, the ``strange term" $\mu u^0$ still appears in the left-hand side while the source term $F(x,u^0)$ is not modified in the right-hand side.

\end{abstract}

\centerline{CONTENTS}
\begin{enumerate}
\item[$1$.] Introduction  \dotfill \pageref{introduction}
\item[$2$.] Notation and assumptions  \dotfill \pageref{assumptions}
\begin{enumerate}
\item[$2.1$.] Notation  \dotfill \pageref{notation}
\item[$2.2$.] The matrix $A(x)$ and the function $F(x,s)$  \dotfill \pageref{2.1}
\item[$2.3$.]  The perforated domains $\oeps$  \dotfill \pageref{2.2}
\end{enumerate}
\item[$3$.] Definition of a solution to the singular semilinear problem in~$\oeps$  \dotfill \pageref{secdefi}
\begin{enumerate} 
\item[$3.1$.]  The space $\mathcal V(\oeps)$ of test functions and the definition of a solution to the problem in $\oeps$  \dotfill \pageref{secdefiV}
\item[$3.2$.]  Statements of existence, stability and uniqueness results  for the problem in $\oeps$  \dotfill \pageref{stat}
\end{enumerate}
\item[$4$.] Definition of a solution to the homogenized singular semilinear problem\break in $\Omega$  \dotfill \pageref{sub82}
\item[$5$.] Statement of the homogenization result  for the singular semilinear problem in~$\oeps$  \dotfill \pageref{52}
\item[$6$.] Definition of the function $z^\eps$, and the strong convergence of $\chi_{_{\oeps}}$   \dotfill \pageref{section6}
\begin{enumerate}
\item[$6.1$.] Definition of the function $z^\eps$, a variant of the test function $w^\eps$  \dotfill \pageref{6.1}
\item[$6.2$.] Strong convergence of the sequence $\chi_{_{\oeps}}$  in $L^1(\Omega)$  \dotfill \pageref{6.2}
\end{enumerate}
\item[$7$.] A priori estimates for the solutions to the singular semilinear problem\break in $\oeps$  \dotfill \pageref{section7}
\item[$8$.] Proof of the Homogenization Theorem~\ref{homogenization}  \dotfill \pageref{section8}
\end{enumerate}

\section{Introduction}
\label{introduction}
The present paper deals with the homogenization of the following strongly\break singular semilinear problem posed in perforated domains $\oeps$:
\begin{equation}
\tag{{\theequation}$^\eps$}
\label{eqprimai}
\begin{cases}
u^\eps \geq 0 & \mbox{in} \; \oeps,\\
\displaystyle - div \,A(x) D u^\eps  = F(x,u^\eps) & \mbox{in} \; \oeps,\\
u^\eps = 0 & \mbox{on} \; \partial \oeps.\\
\end{cases} 
\end{equation}
Here $A(x)$ is a $N\times N$ bounded coercive matrix, $F(x,s)$ is a Carath\'erodoy function  
$
F(x,s):\Omega\times [0,+\infty[\rightarrow [0,+\infty]
$ 
which possibly has a very strong singularity at $s=0$; an example of such a function $F(x,s)$ is 
\begin{align}
\label{21bisi}
\dys F(x,s)= f(x) \frac{\left(a+ \sin (\frac 1s) \right)}{\exp (- \frac {1}{s})} + g(x) \dys\frac{\left(b+ \sin (\frac 1s) \right)}{s^\gamma} + l(x)\,\,
 \dys \mbox{a.e. } x\in\Omega,\, \forall s>0,
\end{align}
where $\gamma>0,$ $a>1,$ $b>1$ and where the functions $f$, $g$ and $l$ are nonnegative; another example is given in \eqref{21bis} below. The precise assumptions that we actually make on the function $F(x,s)$ are given in Subsection~2.1 below. On the other hand, the open sets $\oeps$ are obtained by removing many small closed holes from a fixed open set $\Omega\subset \mathbb{R}^N$, $N\geq 2$. The model example is the case where a bounded open set $\Omega$ is perforated by small holes which are closed balls of radius $r^\eps$ with 
\begin{equation}
\label{radius}
\begin{cases}
r^\eps=C_0 \eps^{N/(N-2)}\,\mbox{ if } N\geq 3,\\
r^\eps=\exp( -C_0/\eps^2)\, \mbox{ if } N= 2,
\end{cases}
\end{equation}
which are periodically distributed in $\mathbb{R}^N$ at the vertices of an $N$-dimensional lattice of cubes of size $2\eps$.  The general framework that we will use for $\oeps$ in the present paper is (a slight generalization of) the one studied by D.~Cioranescu and F.~Murat in \cite{CM} (see also \cite{MH} and \cite{DG}); it will be described in details in Subsection~\ref{2.2} below.

Note that  in \eqref{eqprimai} the homogeneous Dirichlet boundary condition is imposed on the whole boundary of $\oeps$, which includes the boundary of all the holes. In the classical case where the singular semilinear term $F(x,u^\eps)$ is replaced by a fixed source term $f(x)\in L^2(\Omega)$ which does not depend on $u^\eps$, the homogenization in the framework of \cite{CM} of problem \eqref{eqprimai} leads to a problem where ``a strange term" $\mu u^0$ appears in the left-hand side, where $\mu$ is a bounded nonnegative measure of $H^{-1}(\Omega)$ which depends on the holes and which is actually the asymptotic memory of them.

In \cite{GMM1} we treated the case of problem \eqref{eqprimai} where the singularity at $s=0$ is mild, namely the case where 
\begin{align}
\label{mild}
\dys 0\leq F(x,s)\leq h(x)(\frac{1}{s^{\gamma}}+1)
\dys \mbox{ a.e. } x\in\Omega \mbox{ for some } 0\leq\gamma\leq 1.
\end{align}
In that paper we proved, for $\eps$ fixed, existence, stability and uniqueness results for the solution to \eqref{eqprimai}, as well as the homogenization result for perforated domains of the type described above. In this case where \eqref{mild} is satisfied, the solutions to problem \eqref{eqprimai} belong to $H_0^1(\oeps)$, the equation is intended in the usual weak sense (or more exactly in a slight variant of it), and the test functions that we use belong to $H_0^1(\oeps)$.

In contrast, the purpose of the present paper is to treat the case with strong singularities, namely the case where $\gamma>1$ in \eqref{mild}, or more generally where $F(x,s)$ can exhibit any type of singularity at $s=0$ (see for example \eqref{21bisi} above and example \eqref{21bis} below). In this case the solutions $u^\eps$ to \eqref{eqprimai} do not in general belong to $H_0^1(\oeps)$ (see \cite{LM}), but only to $H_{\mbox{\tiny loc}}^1(\oeps)$, even if $u^\eps$ vanishes in some sense on $\partial\oeps$. This induces significant difficulties in order to define a convenient notion of solution, on the first hand in the description of the space which the solution has to belong to, and on the second hand in the definition of the space $\mathcal{V}(\oeps)$ of test functions to be used in the equation. In \cite{GMM2} we introduced a notion of solution by defining non standard spaces for the solution and for the test functions, and by writing the equation like in the definition of solutions by transposition introduced  by J.-L.~Lions and E.~Magenes and by G.~Stampacchia. This framework, which is recalled in Subsection~3.1 below, allowed us to prove in \cite{GMM2} existence, stability and uniqueness results. The present paper uses this framework and can therefore be considered as a continuation of \cite{GMM2}. It is also a confirmation of the fact that the framework introduced in \cite{GMM2} is robust.

In the present paper we prove that, as in the case studied in \cite{CM} where $F(x,s)$ depends only on $x$, the ``strange term" $\mu u^0$ appears in the left-hand side of the homogenized problem while the source term $F(x,u^0)$ is not modified in the right-hand side. In other terms, see Theorem~\ref{homogenization} below, we prove that a subsequence of solutions $u^\eps$ to \eqref{eqprimai} converges, in a convenient sense, to a  solution $u^0$ to the homogenized problem
\begin{equation}
\label{1555}
\begin{cases}
u^0 \geq 0 & \mbox{in} \; \Omega,\\
\displaystyle - div \,A(x) D u^0+\mu u^0  = F(x,u^0) & \mbox{in} \; \Omega,\\
u^0 = 0 & \mbox{on} \; \partial \Omega.\\
\end{cases} 
\end{equation}
Note that the definition of solution that we use for the solution $u^0$ to \eqref{1555} is a variant of the definition introduced in \cite{GMM2}. This definition is recalled in Section~\ref{sub82} below (see also Section~6 of \cite{GMM4}).

 This homogenization result was not a priori obvious, since the holes ``tend to invade the whole of $\Omega$" (see Remark~\ref{66662} below) and since the source term $F(x,u^\eps)$ has a singular behaviour at the boundary of the holes. 

The method of the proof consists in merging the methods of \cite{CM} and of \cite{GMM2}. This however presents some difficulty, since the solution $u^\eps$ to \eqref{eqprimai} in general does not belong to $H_0^1(\oeps)$. This leads us (see Section~\ref{section6} below) to modify the test function $w^\eps$ used in \cite{CM} (which is more or less the difference between $1$ and the capacitary potential of the holes in $\Omega$) by introducing a variant $z^\eps$ of it, which now belongs to the space $\mathcal{V}(\oeps)$ of test functions introduced in \cite{GMM2}.  

In the best of our knowledge, there are only a very few papers concerned with homogenization in the context of this type of singular semilinear problems. In the paper \cite{BC} the authors deal with the case where, in a fixed domain $\Omega$, the matrices $A^\eps(x)$ wildly vary with $\eps$, remaining uniformly bounded and coercive. For that they use the framework introduced in \cite{BO} which is based on the use of strong maximum principle and on the assumption that the function $F(x,s)$ is nonincreasing in $s$. Note that these properties are never used in the present paper, and neither in \cite{GMM1}, \cite{GMM2}, \cite{GMM3}  and \cite{GMM4}. On the other hand, in the paper \cite{DoG} the authors consider the homogenization of singular semilinear problems posed in a domain divided in two parts separated by an oscillating interface. Lastly, in the paper \cite{GVV}, the authors study the homogenization in infinite cylinders perforated with small holes with Dirichlet boundary condition. In contrast, there are many papers concerned with existence and uniqueness of solutions to \eqref{eqprimai} for $\eps$ fixed. Let us just quote, inter alia, \cite{BO}, \cite{CRT}, \cite{LM}, \cite{OP} and \cite{St}. 

To conclude this Introduction, let us mention that in the present context of strongly singular semilinear problems we have not been able to prove a corrector result, while we were able to do it in \cite{GMM1} in the context of mild singularities. The corrector result thus remains for us an open problem in the case of strong singularities.

The plan of the present paper is as follows: In Section~\ref{assumptions} we give the assumptions that we make on the matrix $A(x)$, on the function $F(x,s)$ and on the sequence of perforated domains $\oeps$. In Section~\ref{secdefi} we recall the definition introduced in \cite{GMM2} of the solution to the strongly singular semilinear problem posed in $\oeps$, and the results of existence, stability and uniqueness obtained in \cite{GMM2}. Note that these solutions satisfy a priori estimates which are recalled in Section~\ref{section7}. In Section~\ref{sub82}, we recall the definition given in \cite{GMM4} of the solution to the homogenized problem with a strange term \eqref{1555} (this definition is a variant of the definition given in \cite{GMM2}). In Section~\ref{52}, we state the main result of the present paper, namely the homogenization result for problem \eqref{eqprimai}. This result is proved in Section~\ref{section8}. An important tool for this proof, namely the function $z^\eps$ which replaces here the function $w^\eps$ used in \cite{CM}, is defined in Section~\ref{section6}.

\section{Assumptions and Notation}
\label{assumptions}

As said in the Introduction, in this paper we deal with the asymptotic behaviour, as $\eps$ tends to zero, of solutions to the singular semilinear elliptic problem
\begin{equation}
\tag{{\theequation}$^\eps$}
\label{eqprimad}
\begin{cases}
u^\eps \geq 0 & \mbox{in} \;\oeps,\\
\displaystyle - div \,A(x) D u^\eps  = F(x,u^\eps) & \mbox{in} \; \oeps,\\
u^\eps = 0 & \mbox{on} \; \partial \oeps,\\
\end{cases} 
\end{equation}
where $F(x,s)$ is possibly singular at $s=0$, where $u^\eps$ satisfies the homogenous Dirichlet boundary condition on the whole of the boundary of $\oeps$, and where  $\Omega^\eps$ is a perforated domain obtained by  removing many small holes from a given open bounded set  $\Omega$ in $\mathbb{R}^N$, with a repartition of those many small holes producing a ``strange term" when $\eps$ tends to $0$.

After the brief Subsection~{2.1} dealing with some notation, we begin by giving in Subsection~{2.2} the assumptions on the matrix $A(x)$ and on the function $F(x,s)$; then in Subsection~{2.3}  we describe the geometry of the perforated domains and (a slightly generalization of) the framework introduced in \cite{CM} for treating this problem when the right-hand side is $F(x,u)=f(x)$ in $L^2(\Omega)$. 
\vskip1cm

\noindent {\bf 2.1. Notation}
\label{notation}
$\mbox{}$

In this paper $\Omega$ denotes a bounded open subset of $\mathbb{R}^N$. 

We denote by $\mathcal{D}(\Omega)$ the space of the $C^\infty(\Omega)$ functions  whose support is compact and included in $\Omega$, and by $\mathcal{D}'(\Omega)$ the space of distributions on $\Omega$.

We denote by $\mathcal{M}_b^+(\Omega)$ the space of nonnegative bounded Radon measures on $\Omega$.


Since $\Omega$ is bounded, $\|D w\|_{L^2(\Omega)^N}$ is a norm which is equivalent to $\|w\|_{H^1(\Omega)}$ on $H_0^1(\Omega)$. We set 
$$
\|w\|_{H_0^1(\Omega)}=\|D w\|_{(L^2(\Omega))^N}\, \,\,\forall w\in H_0^1(\Omega).
$$

For every $s\in \mathbb{R}$ and every $k>0$ we define as usual 
$$s^+=\max\{s,0\},\,\,s^-=\max\{0,-s\},$$
$$T_k(s)=\max\{-k,\min\{s,k\}\},\,\,\,\,\,G_k (s)= s- T_k (s).$$

For any measurable function $l:x\in\Omega\rightarrow l(x)\in [0,+\infty]$ we denote
$$
\{l=0\}=\{x\in\Omega: l(x)=0\},\,\,\,\,\,\{l>0\}=\{x\in\Omega: l(x)>0\}.
$$

Finally, in the present paper, we denote by $\varphi$ functions which belong to\break $H_0^1(\Omega)\cap L^\infty(\Omega)$, 
while we denote by $\phi$ functions which belong to~$\mathcal{D}(\Omega)$.

\vskip1cm

\noindent {\bf 2.2. The matrix $A(x)$ and the function $F(x,s)$}
\label{2.1}
$\mbox{}$

In this Subsection, we give the precise assumptions that we make on the data of problem \eqref{eqprimad}.

We assume that 
\begin{equation}
\label{Omega}
\Omega\quad  \mbox{is an open bounded set of} \quad  \mathbb{R}^N,\, N\geq2,
\end{equation}
 (no regularity is assumed on the boundary $\partial\Omega$ of $\Omega$), that the matrix $A$ is bounded and coercive, i.e. satisfies
\begin{equation}\label{eq0.0}
A(x)\in (L^\infty(\Omega))^{N\times N},\,\,
\exists \alpha>0,\, \,A(x)\geq \alpha I \,\,\,\,\,{\rm a.e.}\,\,  x\in\Omega,
\end{equation}
and that the function $F$ satisfies
\begin{equation}
\label{car} 
\begin{cases}
F: (x,s)\in\Omega\times [0, +\infty[ \rightarrow F(x,s)\in [0, +\infty] \,\, \text {is a Carath\'eodory function},\\ \mbox{i.e. } F \mbox{ satisfies}\\
i)\,\forall s\in [0,+\infty[,\, x\in\Omega\to F(x,s)\in [0,+\infty] \mbox{ is measurable},\\
ii)\, \mbox{for a.e. } x\in\Omega, \, s\in [0,+\infty[\rightarrow F(x,s)\in [0,+\infty] \mbox{ is continuous},
\end{cases}
\end{equation}
\begin{equation}
\label{eq0.1}
\begin{cases}
i) \,\exists h,  h(x)\geq 0 \, \, {\rm a.e.}\,\, x \in \Omega, h \in L^r(\Omega),\\
\dys\mbox{with }\, r=\frac{2N}{N+2}\,\,\,\text {if}\,\, N\geq 3,\,\, r>1 \,\,\text{if} \,\, N=2,\\
ii)\,\exists\Gamma: s\in[0,+\infty[ \rightarrow \Gamma(s)\in [0,+\infty[,  \,\,\Gamma\in C^1([0,+\infty[),\\
\mbox{such that } \Gamma(0)=0 \mbox{ and } \Gamma'(s)>0\,\, \forall s>0,\\
iii) \, \displaystyle 0 \leq F(x,s)\leq  \frac{h (x)}{\Gamma(s)}\,\,  \mbox{a.e. } x \in \Omega, \forall s>0.
\end{cases}
\end{equation}

\begin{remark}
\begin{itemize}[leftmargin=*]
\mbox{ }
\item {\it i}) Note that in the whole of the present paper we assume that
$$
N\geq2,
$$  
(see Remark~\ref{remdim} below).
\smallskip

\item {\it ii}) Note that the matrix $A(x)$ and the function $F(x,s)$ are defined for $x\in\Omega$ and not only for $x\in\oeps$.

\smallskip

\item {\it iii}) The function $F(x,s)$ can have a very wild behaviour in $s$ when $s$ tends to zero. A possible example is given by \eqref{21bisi} above, or more generally by 
\begin{align}
\label{21bis}
\dys F(x,s)= f(x) \frac{\left(a+ \sin (S(s)) \right)}{\exp (- S(s))} + g(x) \dys\frac{\left(b+ \sin (\frac 1s) \right)}{s^\gamma} + l(x)\,\,
 \dys \mbox{a.e. } x\in\Omega,\, \forall s>0,
\end{align}
where $\gamma>0$ $a>1$, $b>1$, 
where the function $S$ satisfies 
\begin{equation}
\label{eq.2.2bis} 
S \in C^1(]0,+\infty[),     \quad   S'(s) < 0 \,\, \forall s>0,   \quad   S(s) \to + \infty  \,\, \mbox{as }  s \to 0,
\end{equation}
 and where the functions $f,g$ and $l$ are nonnegative and belong to $L^r(\Omega)$ with $r$ defined by \eqref{eq0.1} above (see Remark~2.1 {\it viii}) of \cite{GMM2}).

\smallskip

\item {\it iv}) The function $F=F(x,s)$ is a nonnegative Carath\'eodory function with values in $[0,+\infty]$ and not only in $[0,+\infty[$. But, in view of conditions (\ref{eq0.1} {\it ii}) and\break (\ref{eq0.1} {\it iii}), for almost every $x\in\Omega$, the function $F(x,s)$ can take the value $+\infty$ only when $s=0$ (or, in other terms, $F(x,s)$ is finite for almost every $x\in\Omega$ when $s>0$).

\smallskip
\item {\it v}) Note that the growth condition (\ref{eq0.1} {\it iii}) is stated for every $s>0$, while in \eqref{car}  $F$ is supposed to be a Carath\'eodory function defined for $s$ in $[0,+ \infty[$ and not only in $]0,+\infty[$. Indeed an indeterminacy $\dys\frac 00$ appears in $\dys \frac{h(x)}{\Gamma(s)}$ when $h(x)=0$ and $s=0$, while the growth and Carath\'eodory assumptions \eqref{eq0.1} and \eqref{car} imply that 
$$
F(x,s)=0 \quad \forall s\geq 0 \quad \mbox{a.e. on } \{x\in\Omega:h(x)=0\}.
$$
On the other hand, when $h$ is assumed to satisfy  $h(x)>0$ for almost every $x\in\Omega$, one can write  (\ref{eq0.1} {\it iii}) for every $s\geq 0$.
\smallskip 
\item {\it vi}) The function $h$ which appears in hypothesis (\ref{eq0.1} {\it i}) is an element of $H^{-1}(\Omega)$. Indeed, when $N\geq 3$, the exponent $ r=\frac{2N}{N+2}$ is nothing but the H\"older's conjugate $(2^*)'$ of the Sobolev's exponent $2^*$, i.e. 
\begin{equation}
\label{2333}
\mbox{when } N\geq 3,\quad\frac 1r=1-\frac{1}{2^*}, \mbox{ where } \frac{1}{2^*}=\frac 12-\frac 1N.
\end{equation}
Making an abuse of notation, we will set 
\begin{align}
\label{212bis}
2^*=\mbox{any }p \mbox{ with } 1<p<+\infty \mbox{ when } N=2.
 \end{align}
With this abuse of notation, $h$ belongs to $L^r(\Omega)=L^{(2^*)'}(\Omega)\subset H^{-1}(\Omega)$ for all $N\geq 2$ since $\Omega$ is bounded. This result is indeed a consequence of Sobolev's and Trudinger Moser's inequalities, which (with the above abuse of notation) assert that 
\begin{equation}
\label{2334}
\|v\|_{L^{2^*}(\Omega)}\leq C_S \|Dv\|_{(L^2(\Omega))^N}\quad \forall v\in H_0^1(\Omega) \mbox{ when } N\geq 2,
\end{equation}
where $C_S=C_S(N)$ when $N\geq 3$ and $C_S=C_S(p,\Omega)$ when $N=2$. In the latest case, for $p$ given with $1<p<+\infty$, the constant $C_S=C_S(p,\Omega)$ is bounded independently of $\Omega$ when $\Omega\subset Q$, for $Q$ a bounded open set of $\mathbb{R}^2$.  
\smallskip

\item {\it vii}) In Section~5 of \cite{GMM1} we performed the homogenization of problem \eqref{eqprimad} in the case where $F(x,s)$ has a mild singularity at $s=0$, namely in the case where in (\ref{eq0.1} {\it iii}) the function $F(x,s)$ satisfies 
\begin{equation}
\label{mildbis}
0\leq F(x,s)\leq h(x)\left(\frac{1}{s^\gamma}+1\right) \quad \mbox{with} \quad 0<\gamma\leq 1.
\end{equation}
This is a particular case of the general case treated in the present paper, but that case is easier to treat since the solution $u^\eps$ to \eqref{eqprimad} belongs to $H_0^1(\oeps)$ when \eqref{mildbis} holds true.
  This property allowed us to prove in that case a corrector result when the matrix $A(x)$ is symmetric and when $u^0$ further belongs to $L^\infty(\Omega)$, see Theorem~5.5 of \cite{GMM1}.
\end{itemize}
\medskip

Many other remarks can be made about the function $F(x,s)$, and we refer the reader to Section~2 of \cite{GMM2} for them.
\qed

\end{remark}

\vskip 1cm

\noindent {\bf 2.3. The perforated domains $\oeps$}
\label{2.2}
$\mbox{}$

In order to obtain the domain $\oeps$, we perforate the fixed domain $\Omega$ (see \eqref{Omega}) in a way that we describe now. According to (a slight generalization of) the setting presented in \cite{CM}, we consider here, for every $\eps$ which takes its values in a sequence of positive numbers which tends to zero, a finite number $n(\eps)$ of closed sets $T_i^\eps$ of $\mathbb{R}^N$,  $1\leq i \leq n(\eps)$, which are the holes. The domain $\Omega^\eps$ is defined by removing these holes $T_i^\eps$ from $\Omega$, that is by setting
\begin{equation}
\label{Omegaeps}
\Omega^\eps=\Omega- \bigcup_{i=1}^{n(\eps)} T_i^\eps.
\end{equation}

Here, as well as everywhere in the present paper, for every function $y^\eps$ in $L^2(\Omega)$, we define $\reallywidetilde{y^\eps}$ as the extension by zero of $y^\eps$ to $\Omega$, namely by
\begin{equation}
\label{defex}
\dys\reallywidetilde{y^\eps}(x)=\begin{cases}
y^\eps(x) & \text{in}\ \Omega^\eps,\\
\dys 0 &\text{in } \bigcup_{i=1}^{n(\eps)}T_i^\eps;
 \end{cases}
\end{equation}
then $\reallywidetilde{y^\eps}\in L^2(\Omega)$ and $\|\reallywidetilde{y^\eps}\|_{L^2(\Omega)}=\|y^\eps\|_{L^2(\oeps)}$; moreover  
\begin{equation}
\label{57bis}
\dys\mbox{if } y^\eps\in H_0^1(\oeps),\mbox{ then } \reallywidetilde{y^\eps}\in H_0^1(\Omega)\,\,
\dys\mbox{with } \reallywidetilde{D y^\eps}=D\reallywidetilde{y^\eps}\mbox{ and } \|\reallywidetilde{y^\eps}\|_{H_0^1(\Omega)}=\|y^\eps\|_{H_0^1(\oeps)}.
\end{equation}

 We suppose that the sequence of domains $\oeps$ is such that there exist a sequence of functions $w^\eps$, a distribution $\mu\in\mathcal{D}'(\Omega)$ and two sequences of distributions\break $\mu^\eps\in\mathcal{D}'(\Omega)$ and $\lambda^\eps\in\mathcal{D}'(\Omega)$ such that
\begin{equation} \label{cond1}
w^\eps\in H^1(\Omega)\cap L^\infty(\Omega),
\end{equation}
\begin{equation} \label{cond1bis}
0\leq w^\eps\leq 1\mbox{ a.e. } x\in\Omega,
\end{equation}
\begin{equation}
\label{cond2}
\forall \varphi \in H_0^1(\Omega)\cap L^\infty(\Omega),\,  w^\eps \varphi \in H_0^1(\Omega^\eps) \mbox{ and } w^\eps\varphi= \reallywidetilde{w^\eps\varphi} \mbox{ in } \Omega,
\end{equation}
\begin{align}
\label{cond3}
w^\eps\rightharpoonup 1\text{  in } H^1(\Omega) \mbox{ weakly, in } L^\infty(\Omega) \mbox{ weakly-star and a.e. in } \Omega \mbox{ as } \eps\to 0, 
\end{align}
\begin{equation}
\label{cond4}
\mu\in H^{-1}(\Omega),
\end{equation}
\begin{align}
\label{cond5}
\begin{cases}
\dys-div \, {}^t\!A(x)Dw^\eps=\mu^\eps-\lambda^\eps\mbox{ in } \mathcal{D}'(\Omega),\\
 \mu^\eps\in H^{-1}(\Omega),\, \lambda^\eps\in H^{-1}(\Omega),\\
\mu^\eps\geq 0 \mbox{ in } \mathcal{D}'(\Omega),\\
\dys \mu^\eps \rightarrow \mu\mbox{ in } H^{-1}(\Omega)\mbox{ strongly},\\
\dys\langle \lambda^\eps,\reallywidetilde{y^\eps}\rangle_{H^{-1}(\Omega),H_0^1(\Omega)}=0\,\,\, \forall y^\eps\in H_0^1(\Omega^\eps).
\end{cases}
\end{align}

\bigskip

\noindent {\bf The model example for $\oeps$}

The prototype of the examples where assumptions \eqref{cond1}, \eqref{cond1bis},  \eqref{cond2}, \eqref{cond3}, \eqref{cond4} and \eqref{cond5} are satisfied is the case where the matrix $A(x)$ is the identity (and where therefore the operator is the Laplace's operator $-div\, A(x)D=-\Delta$), where $\Omega\subset\mathbb{R}^N$, $N\geq 2$, and where the holes $T_i^\eps$ are balls of radius $r^\eps$ with $r^\eps$ given by
\begin{equation*}
\begin{cases}
r^\eps=C_0 \eps^{N/(N-2)}\,\mbox{ if } N\geq 3,\\
\eps^2\log r^\eps\rightarrow -C_0\, \mbox{ if } N= 2,
\end{cases}
\end{equation*}
for some $C_0>0$ (taking $r\dys^\eps=\exp (- {C_0}/{\eps^2})$ is the model case for $N=2$) which are periodically distributed at the vertices of an $N$-dimensional lattice of cubes of size $2\eps$; in this case the measure $\mu$ is given by 
\begin{equation*}
\begin{cases}
\mu=\dys\frac{S_{N-1}(N-2)}{2^N}C_0^{N-2}\,\mbox{ if } N\geq 3,\\
\mu=\dys\frac{2\pi}{4}\frac{1}{C_0}\, \mbox{ if } N= 2,
\end{cases}
\end{equation*}
see e.g. \cite{CM} and \cite{MH} for more details, and for other examples, in particular for the case where the holes have a different form and/or are distributed on a manifold.
\qed

\begin{remark}
\label{remdim}
In dimension $N=1$, there is no sequence $w^\eps$ which satisfies \eqref{cond2} and \eqref{cond3} whenever for every $\eps$ there exists at least one hole $T_{i_\eps}^\eps$ with $T_{i_\eps}^\eps\cap \overline{\Omega}\not=\O$, see Remark~5.1 of \cite{GMM1} for more details. This is the reason why we assume in the present paper that $N\geq 2$.
\qed
\end{remark}

\bigskip
\noindent {\bf Some properties of $w^\eps$ and  $\mu$}

One deduces from the second assertion of \eqref{cond2} that\footnote{erratum corrige: please note that in equation (5.4) of \cite{GMM1} we should have added the requirement ``and $w^\eps\psi=\reallywidetilde{w^\eps \psi}$" (as it is done in \eqref{cond2} in the present paper).} 
\begin{equation}
\label{57ter}
w^\eps=0 \,\text{ in } \bigcup_{i=1}^{n(\eps)}T_i^\eps; 
\end{equation}
more precisely, \eqref{cond2} means that for every $\eps$ and every $\varphi\in H_0^1(\Omega)\cap L^\infty(\Omega)$, there exists a sequence $\phi_n$ (which depends on $\eps$ and on $\varphi$) such that 
\begin{equation}
\label{29bis}
\phi_n\in\mathcal{D}(\oeps),\,\, \reallywidetilde{\phi_n}\rightarrow \reallywidetilde{w^\eps \varphi} \mbox{ in } H^1(\Omega).
\end{equation}

 On the other hand, taking any $\phi^\eps\in\mathcal{D}(\oeps)$ as test function in the first statement of \eqref{cond5} implies that
\begin{equation}
\label{216bis}
\displaystyle -div\, {}^t\!A(x) D w^\eps  = \mu^\eps\mbox{ in } \mathcal{D}'(\oeps),
\end{equation}
which means that the distribution $\lambda^\eps$ ``only acts on the holes $T_i^\eps$", $i=1,\cdots,n(\eps)$; this fact is also reflected by  the last assertion of \eqref{cond5}.

Taking $v^\eps=w^\eps \phi$, with $\phi\in\mathcal{D}(\Omega)$, $\phi\geq 0$, as test function in the first statement of \eqref{cond5} we have, thanks to the last assertion of \eqref{cond5}, 
$$
\into \phi \, {}^t\!A(x)Dw^\eps Dw^\eps+\into w^\eps\, {}^t\!A(x)D w^\eps D\phi=\langle\mu^\eps,w^\eps\phi\rangle_{H^{-1}(\Omega),H_0^1(\Omega)},
$$
from which using \eqref{cond3} and the fourth statement of \eqref{cond5} we deduce that
\begin{equation}
\label{220bis}
\into \phi A(x)Dw^\eps Dw^\eps\rightarrow \langle\mu,\phi\rangle_{H^{-1}(\Omega),H_0^1(\Omega)}\,\, \forall\phi\in\mathcal{D}(\Omega),\,\phi\geq0,
\end{equation}
and therefore using the coercivity \eqref{eq0.0} that
\begin{equation*}
\mu\geq 0.
\end{equation*}
The distribution $\mu\in H^{-1}(\Omega)$ is therefore also a nonnegative measure. Moreover, using \eqref{220bis}, \eqref{eq0.0} and \eqref{cond3}, one deduces that
\begin{align*}
\begin{cases}
\forall\phi\in \mathcal{D}(\Omega),\phi\geq 0,\\
\dys\into \phi \, d\mu=\langle\mu,\phi\rangle_{H^{-1}(\Omega),H_0^1(\Omega)}= \lim_\eps \into \phi A(x)Dw^\eps D w^\eps\leq C\|\phi\|_{L^\infty(\Omega)},
\end{cases}
\end{align*}
for a constant $C$ which does not depend on $\phi$. Therefore the measure $\mu$ is a finite Radon measure which satisfies $\dys\into d\mu\leq C<+\infty$, or in other terms 
\begin{equation}
\label{mufinite}
\mu\in\mathcal{M}_b^+(\Omega).
\end{equation}

We will therefore use in the present paper the following (well) known result\footnote{the reader who would not enter in this theory could continue reading the present paper assuming in \eqref{cond4} that $\mu$ is a function of $L^r(\Omega)$ (with $r=\frac{2N}{N+2}$ if $N\geq 3$ and $r>1$ if $N=2$) and not only an element of $H^{-1}(\Omega)$.} (see e.g. \cite{DMu} Section~1 and \cite{DMu2} Subsection~2.2 for more details): 
\bigskip

 if $y\in H_0^1(\Omega)$ and if $\nu\in \mathcal{M}^+_b(\Omega)\cap H^{-1}(\Omega)$, then $y$ (or more exactly its quasi-continuous representative for the $H_0^1(\Omega)$ capacity) satisfies  
\begin{align}
\label{57bis2}
\begin{cases}
\dys\forall \nu \in\mathcal{M}_b^+(\Omega)\cap H^{-1}(\Omega), \,\,\forall y\in H_0^1(\Omega),\\ \mbox{ one has }
\dys y\in L^1(\Omega;d\nu)\,\mbox{ with }\,\langle \nu ,y\rangle_{H^{-1}(\Omega),H_0^1(\Omega)}=\into y\,d\nu\,;
\end{cases}
\end{align}
moreover 
\begin{align}
\label{58bis}
\begin{cases}
\dys\forall \nu\in\mathcal{M}_b^+(\Omega)\cap H^{-1}(\Omega), \,\,\forall y\in H_0^1(\Omega)\cap L^\infty(\Omega),\\ \mbox{ one has }
\dys y\in L^\infty(\Omega;d\nu)\, \mbox{ with }\, \|y\|_{L^\infty(\Omega;d\nu)}=\|y\|_{L^\infty(\Omega)};
\end{cases}
\end{align}
therefore when $y\in H_0^1(\Omega)\cap L^\infty(\Omega)$, then $y$ belongs to $L^1(\Omega;d\nu)\cap L^\infty(\Omega;d\nu)$ and therefore to $L^p(\Omega;d\nu)$ for every $p$, $1\leq p\leq +\infty$.

\bigskip
\noindent {\bf The limit problem for a source term in $L^2(\Omega)$}

When one assumes that the holes $T_i^\eps$, $i=1,\cdots,n(\eps)$, are such that the assumptions \eqref{cond1}, \eqref{cond1bis}, \eqref{cond2}, \eqref{cond3}, \eqref{cond4} and \eqref{cond5} hold true, then (see \cite{CM}, or \cite{MH}, or \cite{DG} for a more general framework) for every $f\in L^2(\Omega)$, the (unique) solution $y^\eps$ to the linear problem
\begin{equation}
\label{521}
\begin{cases}
y^\eps\in H_0^1(\oeps),\\
-div\, A(x)Dy^\eps=f\mbox{ in } \mathcal{D}'(\oeps),
\end{cases}
\end{equation}
satisfies 
\begin{equation*}
\reallywidetilde{y^\eps}\rightharpoonup y^0\mbox{ in } H_0^1(\Omega),
\end{equation*}
where $y^0$ is the (unique) solution to 
\begin{equation*}
\begin{cases}
\dys y^0\in H_0^1(\Omega)\cap L^2(\Omega;d\mu),\\
\dys -div\, A(x)Dy^0+\mu y^0=f \mbox{ in } \mathcal{D}'(\Omega),
\end{cases}
\end{equation*}
or equivalently to
\begin{equation}
\label{523}
\begin{cases}
\dys y^0\in H_0^1(\Omega)\cap L^2(\Omega;d\mu),\\
\dys\into A(x)Dy^0Dz+\into y^0 z\,d\mu=\into fz\,\,\, \forall z\in H_0^1(\Omega)\cap L^2(\Omega;d\mu).
\end{cases}
\end{equation}

 Note that the ``strange term" $\mu u_0$ which appears in the limit equation \eqref{523} is the asymptotic memory of the fact that $\reallywidetilde{y^\eps}$ was zero on the holes.
\bigskip

\section{Definition of a solution\\ to the singular semilinear problem in $\oeps$}
\label{secdefi}

In Subsection~\ref{secdefiV} we first recall the definition of a solution to the singular semilinear problem \eqref{eqprimad} which will be used in the present paper; this definition has been introduced in Section~3 of \cite{GMM2}. Then, in Subsection~3.2, we recall the main properties (existence, uniqueness and stability) of such a solution; we will recall in Section~\ref{section7} below a priori estimate which are satisfied by every such solution. All these properties have been stated and proved in Sections~4, 5, 6, and 7 of \cite{GMM2}.
\vskip1cm

\noindent {\bf 3.1.The space $\mathcal V(\oeps)$ of test functions and the definition of a solution to the problem in $\oeps$}
\label{secdefiV}
$\mbox{}$

  In order to recall the notion of solution to problem \eqref{eqprimad} that we will use in the present paper, we recall the definition of the space $\mathcal V(\oeps)$ of test functions and a notation (see Section~3 of \cite{GMM2}).
\begin{definition}(Definition~3.1 of  \cite{GMM2})
\label{spazio}
The space $\mathcal V(\oeps)$ is the space of the functions $v^\eps$ which satisfy
\begin{equation}
             \label{vs}
            v^\eps\in H_0^1(\oeps)\cap L^\infty(\oeps), 
\end{equation}
\begin{equation}
\label{condv}
\begin{cases}
\exists I^\eps\,\text{finite},\, \exists  \hat{\varphi}_i^\eps, \exists \hat{ g}_i^\eps,i\in I^\eps,\exists  \hat{ f}^\eps,\mbox{ with}\\
 \hat{\varphi}_i^\eps \in H_0^1(\oeps)\cap L^\infty(\oeps), \hat{ g}_i^\eps\in (L^2(\oeps))^N , \hat{ f}^\eps\in L^1(\oeps),\\ \mbox{such that } \displaystyle-div\, {}^t\!A(x)Dv^\eps=\sum_{i \in I^\eps} \hat{\varphi}_i^\eps (-div \,\hat{ g}_i^\eps)+\hat{ f}^\eps \,\, \mbox{in } \mathcal{D}'(\oeps).
 \end{cases}
 \end{equation}
 \qed
\end{definition} 

 In the definition of $\mathcal{V}(\oeps)$ we use the notation $\hat{\varphi_i}^\eps$, $\hat{g_i}^\eps$ and $\hat{f}^\eps$ to help the reader to identify the functions which enter in the definition of the functions of $\mathcal{V}(\oeps)$.

Observe that $\mathcal{V}(\oeps)$ is a vector space.

\begin{definition}(Definition~3.2 of  \cite{GMM2})
\label{def32}
 When  $v\in \mathcal V(\oeps)$ with $$\displaystyle-div\, {}^t\!A(x)Dv^\eps=\sum_{i \in I^\eps} \hat{\varphi}_i^\eps (-div \,\hat{ g}_i^\eps)+\hat{ f}^\eps\quad\mbox{in } \mathcal{D}'(\oeps),$$ 
 where $I^\eps$, $\hat{\varphi}_i^\eps$, $\hat{g}_i^\eps$ and $\hat{f}^\eps$ are as in \eqref{condv}, and when $y^\eps$ satisfies $$y^\eps\in H^1_{loc}(\oeps)\cap L^\infty(\oeps) \mbox{ with } \varphi^\eps y^\eps\in H^1_0(\oeps)\, \forall\varphi^\eps\in H^1_0(\oeps)\cap L^\infty(\oeps), $$ we use the following notation: 
\begin{equation}\label{dc}
\langle\langle -div\, {}^t\!A(x)Dv^\eps, y^\eps \rangle\rangle_{\oeps}=  \sum_{i\in I^\eps} \displaystyle \int_{\oeps} \hat{ g}_i^\eps D(\hat{ \varphi}_i^\eps y^\eps)+\displaystyle \int_{\oeps} \hat{ f}^\eps  y^\eps \, \\.
\end{equation}
\qed
\end{definition}

  In notation \eqref{dc}, the right-hand side is correctly defined since\break $\hat{\varphi}_i^\eps y^\eps \in H^1_0(\oeps) $ and since $y^\eps\in L^\infty(\oeps)$. In contrast the left-hand side\break $\langle\langle -div\, {}^t\!A\,Dv^\eps, y^\eps \rangle\rangle_{\oeps}$ is just a notation.
  
   \begin{remark}
 \label{rem34}
 In this Remark we recall some observations which are detailed in Remarks~3.4 and 3.5 of \cite{GMM2}.
 
 \begin{itemize}[leftmargin=*]
\item{$i)$} If $\overline{y^\eps}\in H^1_0(\oeps)\cap L^\infty(\oeps) $, then  $\varphi^\eps \overline{y^\eps}\in H^1_0(\oeps)$ for every $\varphi^\eps \in H^1_0(\oeps)\cap L^\infty(\oeps)$, so that for every $v^\eps\in\mathcal V(\oeps)$, $\langle\langle -div\, {}^t\!A\,Dv^\eps, y^\eps \rangle\rangle_{\oeps}$ is defined. In this case one has  
\begin{equation}\label{classic}
\begin{cases}
\, \forall v^\eps\in\mathcal{V}(\oeps), \, \forall \overline{y^\eps} \in H_0^1(\oeps)\cap L^\infty(\oeps),\\
\langle\langle -div\, {}^t\!A(x)Dv^\eps, \overline{y^\eps} \rangle\rangle_{\oeps}=\langle -div\, {}^t\!A(x)Dv^\eps, \overline{y^\eps}\rangle_{H^{-1}(\oeps),H^1_0(\oeps) }. 
\end{cases}
\end{equation}

\smallskip
\item{$ii)$} If $\varphi^\eps \in H^1_0(\oeps)\cap L^\infty(\oeps)$, then $(\varphi^\eps)^2 \in \mathcal V(\oeps)$, with
\begin{equation}
\label{condv5}
-div\, {}^t\!A(x)D(\varphi^\eps)^2=\hat\varphi^\eps (-div\, \hat{g}^\eps)+\hat{f}^\eps\,\, {\rm in }\,\mathcal{D}'(\oeps),
\end{equation}
with $\hat\varphi^\eps=2\varphi^\eps \in H^1_0(\oeps)\cap L^\infty(\oeps)$, $\hat{g}^\eps=\, {}^t\!A(x)D\varphi^\eps\in (L^2(\oeps))^N$ and\break $\hat{f}^\eps=-2\, {}^t\!A(x)D\varphi^\eps D\varphi^\eps\in L^1(\oeps)$.

More in general,  if  $\varphi_1^\eps$ and $\varphi_2^\eps$ belong to $H^1_0(\oeps)\cap L^\infty(\oeps)$, then $\varphi_1^\eps \varphi_2^\eps$ belongs to $\mathcal V(\oeps)$.

\smallskip
\item{$iii)$} If $\varphi^\eps\in H_0^1(\oeps)\cap L^\infty(\oeps)$ with supp $\varphi^\eps\subset K^\eps$, $K^\eps$ compact, $K^\eps\subset \oeps$, then $\varphi^\eps\in \mathcal V(\oeps)$, since 
\begin{equation}
\label{39bis}
-div\, {}^t\!A(x)D\varphi^\eps=\overline{\phi}^\eps(-div\, {}^t\!A(x)D\varphi^\eps)\, \mbox{ in } \, \mathcal D'(\oeps),
\end{equation}
for every $\overline{\phi}^\eps\in\mathcal D(\oeps)$, with $\overline{\phi}^\eps= 1$ on $K^\eps$.

In particular every $\phi^\eps\in \mathcal{D}(\oeps)$ belongs to $\mathcal V(\oeps)$.
\end{itemize}
\qed
\end{remark}
\smallskip

We now recall the definition of a solution to problem \eqref{eqprimad} that we will use in the present paper.
\begin{definition}(Definition~3.6 of  \cite{GMM2})\label{sol}
Assume that the matrix $A$ and the function $F$ satisfy \eqref{eq0.0}, \eqref{car} and \eqref{eq0.1}. We say that $u^\eps$ is a solution to problem \eqref{eqprimad} if $u^\eps$ satisfies
\addtocounter{equation}{1}
\begin{equation}
\tag{{\theequation}$^\eps$}
\label{sol1}
\begin{cases}
i)\, u^\eps\in L^2(\oeps)\cap H^1_{\mbox{{\tiny loc}}}(\oeps),\\
ii)\, u^\eps(x)\geq 0 \,\,\mbox{ a.e. } x\in\oeps,\\
iii)\, G_k(u^\eps)\in H_0^1(\oeps)\,\, \,\, \forall k>0,\\
iv)\, \varphi^\eps T_k(u^\eps)\in H^1_0(\oeps)\,\, \,\,\forall k>0,\,\, \forall \varphi^\eps \in H^1_0(\oeps)\cap L^\infty (\oeps),
\end{cases}
\end{equation}
\addtocounter{equation}{1}
\begin{equation}
\label{sol2}
\tag{{\theequation}$^\eps$}
\begin{cases}
\displaystyle\forall v^\eps \in \mathcal V(\oeps),\,v^\eps\geq 0,\\ \dys\mbox{with } -div \, {}^t\!A(x)Dv^\eps=\sum_{i \in I^\eps} \hat{ \varphi^\eps_i} (-div\, \hat{ g^\eps_i})+\hat{ f^\eps} \mbox{ in } \mathcal{D}'(\oeps),\\
\mbox{where } \hat{\varphi^\eps_i}\in H_0^1(\oeps)\cap L^\infty(\oeps), \hat{g^\eps_i}\in (L^2(\oeps))^N, \hat{f^\eps}\in L^1(\oeps),\\
\mbox{one has}\\\vspace{0.1cm}
i) \,\dys  \int_{\oeps} F(x,u^\eps) v^\eps<+\infty,\\
ii) \,
\dys \int_{\oeps}\, {}^t\!A(x)Dv^\eps DG_k(u^\eps)+  \displaystyle\sum_{i\in I^\eps} \int_{\oeps} \hat{g^\eps_i} D(\hat{\varphi^\eps_i} T_k(u^\eps))+\int_{\oeps} \hat{f^\eps} T_k(u^\eps)=\\=\langle -div\, {}^t\!A(x)Dv^\eps, G_k(u^\eps) \rangle_{H^{-1}(\oeps),H_0^1(\oeps)}+\langle\langle -div\, {}^t\!A(x)Dv^\eps, T_k(u^\eps)\rangle\rangle_{\oeps}=\\
\displaystyle =\int_{\oeps} F(x,u^\eps) v^\eps \,\,\forall k>0.
\end{cases}
\end{equation}
\qed
\end{definition}

\begin{remark}
\label{phiDu}
When $u^\eps$ satisfies \eqref{sol1}, one has 
\begin{equation}
\label{eq:phiDu}
\varphi^\eps Du^\eps \in (L^2(\oeps))^N \quad \forall \varphi^\eps \in H_0^1(\oeps)\cap L^\infty(\oeps);
\end{equation}
indeed one writes in $(\mathcal{D}'(\oeps))^N$
\begin{align*}
\begin{cases}
\dys \varphi^\eps Du^\eps = \varphi^\eps DT_k(u^\eps)+\varphi^\eps DG_k(u^\eps)=\\
\dys = D(\varphi^\eps T_k(u^\eps))-T_k(u^\eps)D\varphi^\eps+\varphi^\eps DG_k(u^\eps).
\end{cases}
\end{align*}

\qed
\end{remark}

In Definition~\ref{sol}, the requirement \eqref{sol1} is the "space" (which is not a vectorial space) to which the solution should belong, while  requirement (\ref{sol2} {\it ii}) expresses the partial differential equation of \eqref{eqprimad} in terms of (non standard) test functions, in the spirit of the solutions defined by transposition introduced by J.-L. Lions and E. Magenes and by G. Stampacchia.
\smallskip

Indeed, very formally, we have
\begin{align*}
\begin{cases}
\dys``\langle -div\, {}^t\!A(x)Dv^\eps, G_k(u^\eps) \rangle_{H^{-1}(\oeps),H_0^1(\oeps)}=\int_{\oeps}(-div\, {}^t\!A(x)Dv^\eps) \,G_k(u^\eps)=\\
\dys=\int_{\oeps}v^\eps\, (-divA(x)DG_k(u^\eps))",
\end{cases}
\end{align*}
\begin{align*}
\begin{cases}
\dys ``\langle\langle -div\, {}^t\!A(x)Dv^\eps, T_k(u^\eps) \rangle\rangle_{\oeps}=\int_{\oeps}(-div\, {}^t\!A(x)Dv^\eps)\, T_k(u^\eps)=\\
\dys=\int_{\oeps}v^\eps\,(-divA(x)DT_k(u^\eps)",
\end{cases}
\end{align*}
so that (\ref{sol2} {\it ii}) formally means that
$$
``\int_{\oeps} v^\eps\, (-div\, {}^t\!A(x)Du^\eps)=\int_{\oeps} F(x,u^\eps)v^\eps\," \,\, \forall v^\eps \in \mathcal{V}(\oeps),\,v^\eps\geq 0. 
$$
Since every $v^\eps$ can be written as $v^\eps=(v^\eps)^+-(v^\eps)^-$ with $(v^\eps)^+\geq 0$, $(v^\eps)^-\geq 0$,  one has formally (this is formal since we do not know whether $(v^\eps)^+$ and $(v^\eps)^-$ belong to $\mathcal{V}(\oeps)$ when $v^\eps$ belongs to $\mathcal{V}(\oeps)$)
$$
``-div\,A(x)Du^\eps=F(x,u^\eps)"
$$
which is the second statement of \eqref{eqprimad}.

On the other hand, the third assertion of \eqref{sol1} formally implies (this is formal since in the present paper the boundary $\partial\oeps$ of $\oeps$ is not assumed to be smooth) that for every $k>0$, one has ``$G_k(u^\eps)=0$ on $\partial \oeps$\,", i.e. ``$u^\eps\leq k$ on $\partial \oeps$\,", which formally implies that  ``$u^\eps=0$ on $\partial \oeps$\,", which is the third statement of \eqref{eqprimad}.

For other observations about Definition~\ref{sol}, see Remark~3.7 and Proposition 3.8 of~\cite{GMM2}.
\qed

\vskip1cm

\noindent {\bf 3.2. Statements of existence, stability and uniqueness results for the problem in $\oeps$}
\label{stat}
$\mbox{}$

In this Subsection we recall results of existence, stability and uniqueness of the solution to problem \eqref{eqprimad} in the sense of Definition~\ref{sol}. These results have been stated and proved in \cite{GMM2}.

\begin{theorem}[{\bf Existence}] \rm({Theorem~4.1 of \cite{GMM2})}
  \label{EUS}
Assume that the matrix $A$ and the function $F$ satisfy \eqref{eq0.0}, \eqref{car} and \eqref{eq0.1}. Then there exists at least one solution $u^\eps$ to  problem \eqref{eqprimad} in the sense of Definition~\ref{sol}. 
\end{theorem}
\qed
  \begin{theorem}[{\bf Stability}] \rm{(Theorem~4.2 of \cite{GMM2})}
\label{est}
 Assume that the matrix $A$ satisfies assumption \eqref{eq0.0}. Let $F_n$ be a sequence of functions and $F_\infty$ be a function which all  satisfy assumptions
\eqref{car} and \eqref{eq0.1} for the same $h$ and the same $\Gamma$. Assume moreover that 
\begin{equation}
 \label{num1}
 \mbox{\rm{ a.e. }} x\in \Omega, \, F_n(x,s_n)\rightarrow\! F_\infty(x,s_\infty) \mbox{\rm{ if }} s_n\rightarrow s_{\infty}, s_n\geq 0,s_\infty\geq0.
 \end{equation} 
Let $u^\eps_n$ be any solution to problem \eqref{eqprimad}$_n$ in the sense of Definition~\ref{sol}, where \eqref{eqprimad}$_n$ is the problem \eqref{eqprimad} with $F(x,s)$ replaced by $F_n(x,s)$. 

Then there exists a subsequence, still labelled by $n$, and a function $u_\infty$, which is a solution to problem  \eqref{eqprimad}$_\infty$ in the sense of Definition~\ref{sol}, such that (for $\eps$ fixed)
\begin{align}
\label{num2}
\begin{cases}
\dys u^\eps_n\rightarrow u^\eps_\infty \mbox{\rm{ in }} L^2(\oeps)\mbox{\rm{ strongly, in }} H^1_{\mbox{\tiny loc}}(\oeps) \mbox{\rm{ strongly and a.e. in }}\oeps,\\
\dys  G_k(u^\eps_n)\rightarrow  G_k(u^\eps_\infty) \, \mbox{\rm{ in }} H_0^1(\oeps)\, \mbox{\rm{strongly }} \forall k>0,\\
\dys \varphi^\eps T_k (u^\eps_n)\rightarrow \varphi^\eps T_k (u^\eps_{\infty}) \, \mbox{\rm{ in }} H^1_0(\oeps)\, \mbox{\rm{strongly }} \forall k>0,\, 
\forall \varphi^\eps\in H_0^1(\oeps)\cap L^\infty(\oeps).\\
\end{cases}
\end{align}
\end{theorem}
\qed

 \bigskip
 
  Finally, the following uniqueness result holds true when, further to \eqref{car} and \eqref{eq0.1}, the function $F(x,s)$ is assumed to be nonincreasing with respect to $s$, i.e. to satisfy
  \begin{equation}\label{eq0.2}
\dys F(x,s)\leq F(x,t)\,\,\,{\rm a.e.}\,\,\, x \in \Omega,\,\,\forall s,\forall t,\,0\leq t\leq s. 
\end{equation}

 \begin{theorem}[{\bf Uniqueness}] \rm{(Theorem~4.3 of \cite{GMM2})}
 \label{uniqueness}
 Assume that the matrix $A$ and the function $F$ satisfy \eqref{eq0.0}, \eqref{car} and \eqref{eq0.1}.  Assume moreover that the function $F(x,s)$  satisfies assumption \eqref{eq0.2}. Then the solution to problem \eqref{eqprimad} in the sense of Definition~\ref{sol} is unique.
 
 \end{theorem}
 \qed
%
 
 \bigskip

 When assumptions \eqref{eq0.0}, \eqref{car}, \eqref{eq0.1} as well as \eqref{eq0.2} hold true, Theorems~\ref{EUS}, \ref{est} and \ref{uniqueness} together  assert that problem \eqref{eqprimad} is well posed in the sense of Hadamard in the framework of Definition~\ref{sol}.
 
 In Section~7 below, we will recall a priori estimates which are satisfied by every solution to \eqref{eqprimad} in the sense of Definition~\ref{sol}.

\bigskip

 \section{Definition of a solution\\ to the homogeneized singular semilinear problem in $\Omega$}
\label{sub82}

In this Section we recall the definition of the solution to the problem 
\begin{equation}
\label{4.0}
\begin{cases}
\dys u\geq 0, &  \mbox{in} \; \Omega,\\
\displaystyle - div \,A(x) D u +\mu u  = F(x,u) & \mbox{in} \; \Omega,\\
u = 0 & \mbox{on} \; \partial \Omega,\\
\end{cases} 
\end{equation}
when $\mu$ satisfies 
\begin{equation}
\label{muu}
\mu \in  \mathcal{M}_b^+(\Omega)\cap H^{-1}(\Omega).
\end{equation}
This Definition, which has been introduced in Section~6 of \cite{GMM4}, is an adaptation of Definition~\ref{sol} above. 
\begin{definition}(Definition~6.1 of \cite{GMM4})
\label{solgen}
Assume that the matrix $A$, the function $F$ and the Radon measure $\mu$ satisfy \eqref{eq0.0}, \eqref{car}, \eqref{eq0.1} and \eqref{muu}. We say that $u$ is a solution to problem \eqref{4.0} if $u$ satisfies
\begin{equation}
\label{sol1gen}
\begin{cases}
i)\, u\in L^2(\Omega)\cap H^1_{\mbox{{\tiny loc}}}(\Omega),\\
ii)\, u(x)\geq 0 \,\,\mbox{a.e. } x\in\Omega,\\
iii)\, G_k(u)\in H_0^1(\Omega)\quad \forall k>0,\\
iv)\, \varphi T_k(u)\in H_0^1(\Omega)\,\,\forall k>0,\,\, \forall \varphi \in H^1_0(\Omega)\cap L^\infty (\Omega),
\end{cases}
\end{equation}
\begin{align}
\label{sol2gen}
\begin{cases}
\displaystyle\forall v \in \mathcal V(\Omega),\,\,v\geq 0,\\ \dys\mbox{with } -div \, {}^t\!A(x)Dv=\sum_{i \in I} \hat{ \varphi}_i (-div\, \hat{ g}_i)+\hat{ f} \mbox{ in } \mathcal{D}'(\Omega),\\
\mbox{where } \hat{\varphi_i}\in H_0^1(\Omega)\cap L^\infty(\Omega), \hat{g_i}\in (L^2(\Omega))^N, \hat{f_i}\in L^1(\Omega),\\
\mbox{one has}\\\vspace{0.1cm}
i) \,\dys  \into F(x,u) v<+\infty,\\
ii) \dys \into\, {}^t\!A(x)Dv DG_k(u)+  \displaystyle\sum_{i\in I} \into \hat{g_i} D(\hat{\varphi_i} T_k(u))+
\dys\into \hat{f} T_k(u)\, +\into uv d\mu\dys=\\
\dys =\langle -div\, {}^t\!A(x)Dv, G_k(u) \rangle_{H^{-1}(\Omega),H_0^1(\Omega)}\,+\langle\langle -div\, {}^t\!A(x)Dv, T_k(u)\rangle\rangle_{\Omega}\,
+\dys\into uv\, d\mu=\\
\displaystyle =\into F(x,u)\, v\quad \forall k>0.
\end{cases}
\end{align}
\qed
\end{definition}

Note that the term $\dys\into uv d\mu$ has a meaning, as shown in  the following Remark.

\begin{remark}
In (\ref{sol2gen} {\it ii}) the term $\dys\into uv d\mu$ has a meaning since \eqref{sol1gen} and\break $v\in H_0^1(\Omega)\cap L^\infty(\Omega)$ actually imply that 
\begin{equation}
\label{4107}
uv \in L^1(\Omega;d\mu).
\end{equation}
Indeed one can write
$$
uv=T_k(u)v+G_k(u)v,
$$
where $T_k(u)v$ and $G_k(u)$, which belong to $H_0^1(\Omega)$ by (\ref{sol1gen} {\it iv}) and (\ref{sol1gen} {\it iii}), belong to $L^1(\Omega;d\mu)$ in view of \eqref{57bis2}, while $v$, which belongs to $H_0^1(\Omega)\cap L^\infty(\Omega)$, belongs to $L^\infty(\Omega;d\mu)$ in view of \eqref{58bis}.

Actually one can prove (see Section~6 of \cite{GMM4}) that any function $u$ which is a solution to problem \eqref{4.0} in the sense of Definition~\ref{solgen} satisfies the regularity result 
\begin{equation}
\label{4108}
G_k(u)\in L^2(\Omega;d\mu) \quad \forall k>0.
\end{equation}
\noindent On the other hand, since $T_k(u)$ belongs to $L^1(\Omega;d\mu)$ and satisfies $0\leq T_k(u)\leq k$ and since $\Omega$ is bounded, $T_k(u)$ also belongs to $L^\infty(\Omega; d\mu)$ and therefore to $L^2(\Omega; d\mu)$. Together with \eqref{4108} this implies that any solution to problem \eqref{4.0} in the sense of Defintion~\ref{solgen} actually satisfies the regularity result
\begin{equation}
\label{4109}
u\in L^2(\Omega;d\mu).
\end{equation} 
Since $v\in H_0^1(\Omega)\cap L^\infty(\Omega)$ also belongs to $L^2(\Omega;d\mu)$ in view of \eqref{57bis2} and \eqref{58bis}, this again proves \eqref{4107}.

Note however that this second proof of \eqref{4107} uses the fact that $u$ satisfies \eqref{sol1gen} and \eqref{sol2gen}, while the first proof only uses the fact that $u$ satisfies  \eqref{sol1gen}.

\qed
\end{remark}

\begin{remark}
As mentioned in Section~6 of \cite{GMM4}, one can prove, for solutions to problem \eqref{4.0} in the sense of Definition~\ref{solgen}, results of existence, stability and uniqueness which are similar to the results recalled in Subsection~3.2 above for the solutions to problem \eqref{eqprimad} in the sense of Definition~\ref{sol}. Every solution to problem \eqref{4.0} in the sense of Definition~\ref{solgen} moreover satisfies a priori estimates which are similar to the ones recalled in Section~\ref{section7} above, see Section~6 of \cite{GMM4} for more details. 
\qed
\end{remark}

 \section{Statement of the homogenization result\\  for the singular semilinear problem in~$\oeps$}
\label{52}
 $\mbox{}$
 
 The existence Theorem~\ref{EUS} above asserts that when the matrix $A$ and the function $F$ satisfy assumptions \eqref{eq0.0}, \eqref{car} and \eqref{eq0.1}, then for every given $\eps>0$, the singular semilinear problem \eqref{eqprimad} posed in $\oeps$ has at least a solution $u^\eps$ in the sense of Definition~\ref{sol}; moreover (see Theorem~\ref{uniqueness}) this solution is unique if the function $F(x,s)$ also satisfies assumption \eqref{eq0.2}.
 
 The following result, which is the main result of the present paper, asserts that the homogenization process for the singular semilinear problem \eqref{eqprimad} produces a result which is very similar to the homogenization result \eqref{523} above which holds true for the ``classical" problem \eqref{521} when the source term $f$ belongs to $L^2(\Omega)$.

\begin{theorem}{\rm[{\bf Homogenization}]
\label{homogenization}
Assume that the matrix $A$ and the function $F$ satisfy \eqref{eq0.0}, \eqref{car} and \eqref{eq0.1}. Assume also that the sequence of perforated sets $\Omega^\eps$ is such that \eqref{cond1}, \eqref{cond1bis}, \eqref{cond2}, \eqref{cond3}, \eqref{cond4} and \eqref{cond5} are satisfied. For every $\eps>0$, let $u^\eps$ be any solution to problem \eqref{eqprimad} in the sense of Definition~\ref{sol}, or, in other terms, any function which satisfies \eqref{sol1} and \eqref{sol2}, i.e.
\begin{equation}\label{sol1h}
\begin{cases}
i)\, u^\eps\in  L^2(\oeps)\cap H^1_{loc}(\oeps),\\
ii)\, u^\eps(x)\geq 0 \,\,\mbox{\rm{a.e. }} x\in\oeps,\\
iii)\, G_k(u^\eps)\in H^1_0(\oeps) \,\, \forall k>0,\\
iv)\, \varphi^\eps T_k(u^\eps)\in H^1_0(\oeps)\,\,\forall k>0,\,\, \forall \varphi^\eps \in H^1_0(\oeps)\cap L^\infty (\oeps),\\
\end{cases}
\end{equation}
\begin{equation}
\label{condvh}
\begin{cases}
\forall v^\eps\in \mathcal{V}(\oeps), v^\eps\geq 0,\\
 \mbox{\rm{with }} \displaystyle-div\, {}^t\!A(x)Dv^\eps=\sum_{i \in I^\eps} \hat{\varphi}_i^\eps (-div \,\hat{ g}_i^\eps)+\hat{ f}^\eps \,\, \mbox{\rm{in }} \mathcal{D}'(\oeps),\\
  \mbox{\rm{where }}  \hat{\varphi}_i ^\eps \in H_0^1(\oeps)\cap L^\infty(\oeps), \hat{ g}_i^\eps\in L^2(\oeps)^N , \hat{ f}^\eps\in L^1(\oeps),\\
  \mbox{\rm{one has }}\\\vspace{0.1cm}
  i) \dys\int_{\oeps} F(x,u^\eps)v^\eps<+\infty,\\
  \displaystyle ii)\dys\int_{\oeps} {}^t\!A(x) Dv^\eps DG_k(u^\eps)+\sum_{i\in I^\eps} \int_{\oeps} \hat{g}_i^\eps D(\hat{\varphi}_i^\eps T_k( u^\eps))+\into \hat{f}^\eps T_k(u^\eps)=\\
\displaystyle=\langle -div\, {}^t\!A(x)Dv^\eps,G_k(u^\eps)\rangle_{H^{-1}(\oeps), H_0^1(\oeps)}+\langle\langle -div\, {}^t\!A(x)Dv^\eps, T_k(u^\eps)\rangle\rangle_{\oeps }=\\ \dys=\int_{\oeps} F(x,u^\eps)v^\eps \quad \forall k>0.
 \end{cases}
 \end{equation}
Then there exists a subsequence, still denoted by $\eps$, such that  for $\reallywidetilde{u^\eps}$, the extension by zero of $u^\eps$ to $\Omega$ defined by \eqref{defex}, one  has
\begin{equation}
\label{num2h}
\reallywidetilde{u^\eps}\rightharpoonup u^0 \, \mbox{ in } L^{2}(\Omega)\mbox{ weakly and a.e. in } \Omega,\\
\end{equation}
\begin{equation}
\label{53bis}
G_k(\reallywidetilde{u^\eps})\rightharpoonup G_k(u^0) \, \mbox{ in }  H^1_0(\Omega) \mbox{ weakly } \forall k>0,\\
\end{equation}
\begin{equation}
\label{53ter}
\varphi w^\eps T_k(\reallywidetilde{u^\eps})\rightharpoonup \varphi T_k(u^0) \, \mbox{ in }  H^1_0(\Omega) \mbox{ weakly  } \forall k>0,\, \forall \varphi \in H_0^1(\Omega)\cap L^{\infty}(\Omega),\\
\end{equation}
where $u^0$ satisfies \eqref{sol1gen} and \eqref{sol2gen}, or, in other terms, where the limit $u^0$ is a solution to problem \eqref{4.0} in the sense of Definition~\ref{solgen}, i.e.
\begin{equation}\label{sol1h0}
\begin{cases}
i)\, u^0\in  L^2(\Omega)\cap H^1_{loc}(\Omega),\\
ii)\, u^0(x)\geq 0 \,\,\mbox{a.e. } x\in\Omega,\\
iii)\, G_k(u^0)\in H_0^1(\Omega) \,\, \forall k>0,\\
iv)\, \varphi T_k(u^0)\in H^1_0(\Omega)\,\,\forall k>0,\,\,  \forall \varphi \in H^1_0(\Omega)\cap L^\infty (\Omega),\\
\end{cases}
\end{equation}
\begin{align}
\label{condv3h}
\begin{cases}
\forall v\in\mathcal{V}(\Omega),\,\,v\geq0,\\
 \mbox{\rm{with }} \displaystyle-div\, {}^t\!A(x)Dv=\sum_{i \in I} \hat{\varphi}_i (-div \,\hat{ g}_i)+\hat{ f} \,\, \mbox{\rm{in }} \mathcal{D}'(\Omega),\\
\mbox{\rm{where }}  \hat{\varphi}_i \in H_0^1(\Omega)\cap L^\infty(\Omega), \hat{ g}_i\in L^2(\Omega)^N , \hat{ f}\in L^1(\Omega),\\
 \mbox{\rm{one has }}\\\vspace{0.1cm}
 \dys i)\, \dys\into F(x,u^0)v<+\infty,\\
ii)\, \dys\into {}^t\!A(x) DvDG_k(u^0)+\sum_{i\in I} \into \hat{g_i} D(\hat{\varphi_i} T_k(u^0))+\into \hat{f}\,T_k(u^0)\,
\dys+\into u v\,d\mu=\\
\displaystyle=\langle -div\, {}^t\!A(x)Dv,G_k(u^0)\rangle_{H^{-1}(\Omega),H_0^1(\Omega)}+\langle\langle -div\, {}^t\!A(x)Dv, T_k(u^0)\rangle\rangle_{\Omega}
\dys+\into uv\,d\mu=\\\dys= \into F(x,u^0) \,\, \forall k>0.
\end{cases}
\end{align}}
\end{theorem} 

\qed



\begin{remark}
Observe that in the case where assumption \eqref{eq0.2} is made on the function $F(x,s)$, i.e. when $F(x,s)$ is assumed to be nonincreasing with respect to $s$, the solutions $u^\eps$ to problem \eqref{eqprimad} in the sense of Definition~\ref{sol} and $u^0$ to problem \eqref{4.0} in the sense of Definition~\ref{solgen} are unique. In this case there is no need to extract a subsequence in Theorem~\ref{homogenization} and the convergences \eqref{num2h}, \eqref{53bis} and \eqref{condvh} hold true  for the whole sequence $\eps$.
\qed
\end{remark}
\bigskip

\section{Definition of the function $z^\eps$,\\ and the strong convergence of $\chi_{_{\oeps}}$}
\label{section6}
\bigskip
\noindent {\bf 6.1. Definition of the function $z^\eps$, a variant of the test function $w^\eps$}
\label{6.1}
$\mbox{}$
\bigskip

The idea of the proof of the Homogenization Theorem~\ref{homogenization} of the present paper is to combine the ideas of the proof of the Existence Theorem~4.1 of \cite{GMM2} with the ideas of the proof of the Homogenization Theorem~1.2 of \cite{CM}. 
In the latest paper a key tool is the use of the test function $w^\eps \phi$, where $\phi\in \mathcal{D}(\Omega)$ and where $w^\eps$ is defined is in \eqref{cond1}, \eqref{cond1bis},  \eqref{cond2}, \eqref{cond3}, \eqref{cond4} and \eqref{cond5}. Unfortunately, this function does not (seem to) belong to $\mathcal{V}(\oeps)$: 
indeed, the function $w^\eps \phi$ belongs to $H_0^1(\oeps)\cap L^\infty(\oeps)$, but the computation in $\mathcal{D}'(\oeps)$ of $\displaystyle-div\, {}^t\!A(x)D( w^\eps\phi)$ produces four terms, where three of them are in the form required for $w^\eps\phi$ to belong to $\mathcal{V}(\oeps)$, but where the fourth term 
$$
\phi(-div\, {}^t\!ADw^\eps)=\phi(\mu^\eps-\lambda^\eps)=\phi \mu^\eps \mbox{ in } H^{-1}(\oeps),
$$
 is in the form $\phi (-div\, G^\eps)$ for some $G^\eps\in (L^2(\oeps))^N$ (see \eqref{10112bis} below), but not in the requested form $\hat{\varphi}_i^\eps(-div\, \hat{G^\eps})$ with $\hat{\varphi}_i^\eps\in H_0^1(\oeps)\cap L^\infty(\oeps)$, since $\phi$ belongs to $H_0^1(\Omega)\cap L^\infty(\Omega)$ but not to $H_0^1(\oeps)\cap L^\infty(\oeps)$. 
For this reason we introduce in this Section the function $z^\eps$, which is a variant of $w^\eps$ but which is such  that $z^\eps v$ belongs to $\mathcal{V}(\oeps)$ for every $v\in \mathcal{V}(\Omega)$.
Note that the function $z^\eps$ does not (seem to) belong to the smaller (and easier to understand) space of test functions $\mathcal{W}(\oeps)$ introduced in Subsection~4.3 of \cite{GMM4}, which is generated by products $\varphi^\eps \psi^\eps$ with $\varphi^\eps$ and $\psi^\eps$ in $H_0^1(\oeps)\cap L^\infty(\oeps)$. This is actually the reason for which we decided to choose in \cite{GMM2} the framework of the space $\mathcal{V}(\Omega)$ instead of the framework of  the space $\mathcal{W}(\Omega)$.

\begin{proposition}
\label{zeps}
{\rm Assume that \eqref{cond1}, \eqref{cond1bis}, \eqref{cond2}, \eqref{cond3}, \eqref{cond4} and \eqref{cond5} hold true. Then (for a subsequence, as far as the almost everywhere convergence in \eqref{10102} is concerned), there exists a function $z^\eps$ such that
\begin{equation}
\label{10101}
z^\eps\in H^1(\Omega)\cap L^\infty(\Omega),
\end{equation} 
\begin{equation}
\label{10104}
z^\eps-w^\eps \in H_0^1(\oeps),
\end{equation} 
\begin{equation}
\label{10103}
0\leq z^\eps(x)\leq 1 \mbox{ \rm{ a.e.} } x\in\Omega,
\end{equation} 
\begin{equation}
\label{10106}
\,z^\eps v\in\mathcal{V}(\oeps)\quad \forall v\in\mathcal{V}(\Omega),
\end{equation}

\begin{align}
\label{10102}
z^\eps \rightharpoonup 1 \mbox{ \rm{in} } H^1(\Omega) \mbox{ \rm{weakly, in }} L^\infty(\Omega) \mbox{ \rm{weakly-star and a.e. in }} \Omega
\mbox{ as } \eps \rightarrow 0,
\end{align}

\begin{equation}
\label{10105}
-div\, {}^t\!A(x)Dz^\eps=w^\eps \mu^\eps\quad \mbox{\rm{in }} \mathcal{D}'(\oeps).
\end{equation} }
\end{proposition}
\qed
\begin{remark}
Note that in view of \eqref{57ter}, assertion \eqref{10103} implies that in particular
\begin{equation}
\label{66bis}
z^\eps=0 \mbox{ in }  \bigcup_{i=1}^{n(\eps)} T_i^\eps.
\end{equation}

\qed
\end{remark}

\begin{remark}
\label{rem62}
As far as \eqref{10106} is concerned, we will actually prove that
\begin{equation}
\label{10111}
z^\eps \varphi \in H_0^1(\oeps)\cap L^\infty(\oeps)\quad \forall \varphi \in H_0^1(\Omega)\cap L^\infty(\Omega),
\end{equation}
and that if $v$ is such that   
\begin{equation}
\label{10112}
\begin{cases}
\dys v\in\mathcal{V}(\Omega)\\ \mbox{with }
\dys-div \, {}^t\!A(x)Dv=\sum_{i \in I} \hat{ \varphi_i} (-div\, \hat{ g}_i)+\hat{ f} \mbox{ in } \mathcal{D}'(\Omega),\\
\dys\mbox{where } \hat{\varphi_i}\in H_0^1(\Omega)\cap L^\infty(\Omega), \hat{g_i}\in L^2(\Omega)^N, \hat{f}\in L^1(\Omega),\end{cases}
\end{equation}
and if $G^\eps$ is a sequence such that 
\begin{align}
\label{10112bis}
\begin{cases}
\mu^\eps=-div\, G^\eps \mbox{ in }\mathcal{D}'(\Omega),\,\, \mu=-div \, G  \mbox{ in } \mathcal{D}'(\Omega),\\
\mbox{with } G^\eps \rightarrow G \mbox{ in } (L^2(\Omega))^N \mbox{ strongly},
\end{cases}
\end{align}
(note that such a sequence exists since $\mu^\eps$ converges strongly in $H^{-1}(\Omega)$ to $\mu$ in view of the fifth assertion of \eqref{cond5}),
one has 
\begin{align}
\label{1021bis}
\begin{cases}
\dys z^\eps v \in\mathcal{V}(\oeps)\\ \mbox{with}
\dys-div \, {}^t\!A(x)D(z^\eps v)=\\
\dys = \sum_{i\in I} z^\eps \hat{\varphi_i}\,(-div\,\hat{g_i})+w^\eps v\,(-div\, G^\eps)\,
\dys + z^\eps \hat{f}-\, {}^t\!A(x)DvD z^\eps-\, {}^t\!A(x)D z^\eps Dv\\ \mbox{in } \mathcal{D}'(\oeps).
\end{cases}
\end{align}
\qed
\end{remark}
\begin{proof}[\rm{\bf{Proof of Proposition~\ref{zeps}}}]
\mbox{ }

\noindent{\bf First step}. Since $w^\eps\in H^1(\Omega)$ and $\mu^\eps \in H^{-1}(\Omega)$,  the product $w^\eps\mu^\eps$ is, as usual, the distribution on $\Omega$ defined, for every $\phi\in \mathcal{D}(\Omega)$, as
\begin{equation}
\label{1021}
\langle w^\eps \mu^\eps,\phi\rangle_{\mathcal{D}'(\Omega),\mathcal{D}(\Omega)}=\langle \mu^\eps,w^\eps\phi\rangle_{H^{-1}(\Omega),H_0^1(\Omega)}   .
\end{equation}

We claim that actually
\begin{equation}
\label{1022}
w^\eps \mu^\eps \in H^{-1}(\Omega);
\end{equation}
indeed, since $\mu^\eps\geq 0$ in $\mathcal{D}'(\Omega)$ (see the third assertion of \eqref{cond5}), $\mu^\eps$ is a nonnegative Radon measure on $\Omega$, and therefore $\mu^\eps$ belongs to $\mathcal{M}^+_b(\omega)$ for every open set $\omega$ with $\overline{\omega}\subset \Omega$. Taking, for any given $\phi\in\mathcal{D}(\Omega)$, an open set $\omega$ with $\supp\, \phi\subset \omega\subset \overline{\omega}\subset\Omega$, we have,  using \eqref{1021} in $\Omega$ and then \eqref{57bis2} in $\omega$, 
\begin{equation*}
\begin{cases}
\dys\langle w^\eps \mu^\eps,\phi\rangle_{\mathcal{D}'(\Omega),\mathcal{D}(\Omega)}=\langle \mu^\eps,w^\eps\phi\rangle_{H^{-1}(\Omega), H_0^1(\Omega)}=\\
\dys=\langle \mu^\eps,w^\eps\phi\rangle_{H^{-1}(\omega), H_0^1(\omega)}=\int_{\omega} w^\eps \phi\, d\mu^\eps,
\end{cases}
\end{equation*}
and then, using \eqref{cond1bis} and \eqref{57bis2} in $\Omega$, 
\begin{align*}
\begin{cases}
\dys |\langle w^\eps\mu^\eps,\phi\rangle_{\mathcal{D}'(\Omega),\mathcal{D}(\Omega)}|= \left|\int_{\omega} w^\eps \phi\, d\mu^\eps\right|\leq\int_{\omega} w^\eps |\phi|d\mu^\eps\leq \int_{\omega} |\phi| \,d\mu^\eps=\\
\dys=\int_{\Omega} |\phi| \,d\mu^\eps=\dys \langle \mu^\eps,|\phi| \rangle_{H^{-1}(\Omega), H_0^1(\Omega)}
\dys\leq  \|\mu^\eps\|_{H^{-1}(\Omega)} \|\phi\|_{H_0^1(\Omega)} \,\,\forall \phi \in \mathcal{D}(\Omega).
\end{cases}
\end{align*}
This implies that \eqref{1022} holds true with
\begin{equation}
\label{10.7bis}
\|w^\eps \mu^\eps\|_{H^{-1}(\Omega)}\leq \|\mu^\eps\|_{H^{-1}(\Omega)}.
\end{equation}

\bigskip
\noindent{\bf Second step}. Since $w^\eps \mu^\eps$ belongs to $H^{-1}(\Omega)$ by \eqref{1022}, one has
$$
w^\eps \mu^\eps \in H^{-1}(\Omega)\subset H^{-1}(\oeps).
$$
  Applying Lax-Milgram's Lemma then implies the existence (and the uniqueness) of the solution $y^\eps$ to
 \begin{equation}
\label{1020}
\begin{cases}
y^\eps\in H^1(\oeps),\\
y^\eps-w^\eps\in H_0^1(\oeps),\\
-div\, {}^t\!A(x)Dy^\eps=w^\eps \mu^\eps\quad \mbox{in } \mathcal{D}'(\oeps).
\end{cases}
\end{equation}

We now define $z^\eps$ by 
\begin{equation}
\label{defzetaeps}
z^\eps=\reallywidetilde{y^\eps},
\end{equation}
where $\reallywidetilde{y^\eps}$ is the extension by zero of $y^\eps$ to $\Omega$ defined by \eqref{defex}. Then $z^\eps\in H^1(\Omega)$ and $z^\eps$ satisfies \eqref{10104} and \eqref{10105}. 

\bigskip
\noindent{\bf Third step}. We now prove that
\begin{equation}
\label{1026}
0\leq z^\eps(x)\leq w^\eps(x) \mbox{ a.e. } x\in\Omega,
\end{equation}
a fact which in particular implies \eqref{10103} in view of \eqref{cond1bis}, and which completes the proof of \eqref{10101}.

Since one deduces \eqref{66bis} from 
$$
z^\eps=0 \mbox{ and } w^\eps=0 \mbox{ in }  \bigcup_{i=1}^{n(\eps)} T_i^\eps,
$$
(see \eqref{defzetaeps}, \eqref{defex} and \eqref{57ter}), we only have to prove that
\begin{equation}
\label{617bis}
0\leq y^\eps(x)\leq w^\eps(x) \,\mbox{ a.e. } x\in\oeps.
\end{equation}

In order to prove \eqref{617bis}, we first observe that $-(y^\eps)^{-}\in H_0^1(\oeps)$: indeed $-(y^\eps)^{-}\in H^1(\oeps)$ in view of \eqref{1020} and one has
$$
-(y^\eps-w^\eps)^{-}\leq -(y^\eps)^{-}\leq 0 \quad \mbox{a.e. in } \oeps,
$$
where the first inequality results from the facts that the function $-s^{-}$ is nondecreasing and that $w^\eps\geq 0$ (see \eqref{cond1bis}); therefore Lemma~A.1 of \cite{GMM2} implies that $-(y^\eps)^{-}\in H_0^1(\oeps)$. Using  $-(y^\eps)^{-}$ as test function in \eqref{1020} we get, in view of \eqref{cond1bis} and of the third assertion of \eqref{cond5}, 
$$
\int_{\oeps} {}^t\!A(x) Dy^\eps D(-(y^\eps)^{-})=\langle w^\eps \mu^\eps, -(y^\eps)^{-}\rangle_{H^{-1}(\oeps),H_0^1(\oeps)} \leq 0,
$$
which implies that
$$
0\leq y^\eps(x) \mbox{ a.e. } x\in\oeps.
$$

On the other hand, since $y^\eps-w^\eps\in H_0^1(\oeps)$ by \eqref{1020}, using\break $(y^\eps-w^\eps)^+\in H_0^1(\oeps)$ as test function in \eqref{1020} and \eqref{216bis} we get
\begin{align*}
\dys\int_{\oeps} {}^t\!A(x) D(y^\eps-w^\eps) D(y^\eps-w^\eps)^+=
 \langle w^\eps\mu^\eps-\mu^\eps,(y^\eps-w^\eps)^+\rangle_{H^{-1}(\oeps),H_0^1(\oeps)}.
\end{align*}
Since in view of \eqref{cond1bis} and of the third assertion of \eqref{cond5} one has
$$
\langle (w^\eps-1)\mu^\eps,(y^\eps-w^\eps)^+\rangle_{H^{-1}(\oeps),H_0^1(\oeps)}\leq 0,
$$
this implies that 
$$
y^\eps(x)-w^\eps(x)\leq 0 \mbox{ a.e. } x\in\oeps.
$$

We have proved that \eqref{617bis} (and therefore \eqref{1026}) holds true.

\bigskip
\noindent{\bf Fourth step}.  Let us now prove that
\begin{equation}
\label{1015bis}
z^\eps-w^\eps\rightarrow 0 \mbox{ in } H_0^1(\Omega) \mbox{ strongly}.
\end{equation}
Combined with \eqref{cond3} and \eqref{10103}, this will imply \eqref{10102} (for a subsequence, as far as the almost everywhere convergence is concerned).

Using $y^\eps-w^\eps\in H_0^1(\oeps)$ as test function in \eqref{1020} and \eqref{216bis} we get
\begin{align*}
\dys\int_{\oeps} {}^t\!A(x) D(y^\eps-w^\eps) D(y^\eps-w^\eps)
= \langle w^\eps\mu^\eps-\mu^\eps,y^\eps-w^\eps\rangle_{H^{-1}(\oeps),H_0^1(\oeps)}.
\end{align*}

Using the coercivity of the matrix $A$, this implies that
\begin{align}
\label{1011}
\begin{cases}
\alpha\|y^\eps-w^\eps\|_{H_0^1(\oeps)}^2
 \leq\langle w^\eps\mu^\eps-\mu^\eps,y^\eps-w^\eps\rangle_{H^{-1}(\oeps),H_0^1(\oeps)}\leq\\
\leq \left(\|w^\eps \mu^\eps\|_{H^{-1}(\oeps)}+\|\mu^\eps\|_{H^{-1}(\oeps)}\right) \|y^\eps-w^\eps\|_{H_0^1(\oeps)},
\end{cases}
\end{align}
which in view of \eqref{10.7bis}, of the fourth assertion of \eqref{cond5} and of \eqref{57bis} implies that $\|y^\eps-w^\eps\|_{H_0^1(\oeps)}=\|z^\eps-w^\eps\|_{H_0^1(\Omega)}$ is bounded. But $z^\eps-w^\eps$ is also bounded in $L^\infty(\Omega)$ in view of \eqref{10103} and \eqref{cond1bis}. Therefore $(w^\eps-1)(z^\eps-w^\eps)$ is bounded in $H_0^1(\Omega)\cap L^\infty(\Omega)$, and in view of \eqref{cond3}
\begin{equation}
\label{1016bis}
(w^\eps-1)(z^\eps-w^\eps)\rightharpoonup 0 \quad \mbox{in } H_0^1(\Omega) \mbox{ weakly}.
\end{equation}
Writing
\begin{align*}
\begin{cases}
\langle w^\eps \mu^\eps-\mu^\eps, y^\eps-w^\eps\rangle_{H^{-1}(\oeps),H_0^1(\oeps)}=\langle w^\eps\mu^\eps-\mu^\eps,z^\eps-w^\eps\rangle_{H^{-1}(\Omega),H_0^1(\Omega)}=\\
= \langle \mu^\eps, (w^\eps-1)(z^\eps-w^\eps)\rangle_{H^{-1}(\Omega),H_0^1(\Omega)},
\end{cases}
\end{align*}
and using the fact that $\mu^\eps$ tends to $\mu$ in $H^{-1}(\Omega)$ strongly by the fourth assertion of \eqref{cond5}, we deduce \eqref{1015bis} from the first line of \eqref{1011}.

\bigskip
\noindent{\bf Fifth step}. At this point, we have proved the existence of a sequence $z^\eps$ which satisfies \eqref{10101}, \eqref{10104}, \eqref{10103}, \eqref{10102} and \eqref{10105}. Let us now prove that $z^\eps$ satisfies \eqref{10106},  or more precisely \eqref{10111} and \eqref{1021bis} when $v\in\mathcal{V}(\Omega)$ satisfies \eqref{10112}.

Assertion \eqref{10111} follows from the equality
$$
z^\eps\varphi =(z^\eps-w^\eps)\varphi+w^\eps\varphi,
$$
and from the facts that  when $\varphi\in H_0^1(\Omega)\cap L^\infty(\Omega)$, then both $(z^\eps-w^\eps)\varphi$ and $w^\eps \varphi$ belong to $H_0^1(\oeps)\cap L^\infty(\oeps)$ (see \eqref{10104}, \eqref{10103}, \eqref{cond1}, \eqref{cond1bis} and \eqref{cond2}).

On the other hand, using \eqref{10112}, \eqref{1020} and \eqref{10112bis}, we have in $\mathcal{D}'(\oeps)$
\begin{align}
\label{10113}
\begin{cases}
\dys-div \, {}^t\!A(x)D(z^\eps v)\dys=-div(z^\eps \, {}^t\!A(x)Dv)- div(v \, {}^t\!A(x)Dz^\eps)=\\
\dys= z^\eps (-div \, {}^t\!A(x)D v)-\, {}^t\!A(x)DvD z^\eps+\\
\dys +\,v\, (-div \, {}^t\!A(x)D z^\eps)- \, {}^t\!A(x)D z^\eps Dv=\\
\dys = \sum_{i\in I} z^\eps \hat{\varphi_i}(-div\,\hat{g_i})+z^\eps \hat{f}
\dys -\, {}^t\!A(x)DvD z^\eps+\\
\dys + v w^\eps  (-div\, G^\eps) -\, {}^t\!A(x)D z^\eps Dv \mbox{ in } \mathcal{D}'(\oeps),
\end{cases}
\end{align}
which completes the proof of \eqref{1021bis}. This proves \eqref{10106}, in particular since $z^\eps \hat{\varphi_i}$ and $w^\eps v$ belong to  $H_0^1(\oeps)\cap L^\infty(\oeps)$, by \eqref{10111}, \eqref{cond2} and \eqref{cond1}.
\end{proof}
\vskip1cm
\noindent {\bf 6.2. Strong convergence of the sequence $\chi_{_{\oeps}}$ in $L^1(\Omega)$}
\label{6.2}
$\mbox{}$

In this Subsection we prove the following Proposition:

\begin{proposition}{\rm
Assume that the sequence of perforated set $\oeps$ is such that \eqref{cond1}, \eqref{cond1bis}, \eqref{cond2}, \eqref{cond3}, \eqref{cond4} and \eqref{cond5} are satisfied. Then
\begin{equation}
\label{226}
\chi_{_{\oeps}}\rightarrow 1 \mbox{ in } L^1(\Omega) \mbox{ strongly as } \eps \rightarrow 0.
\end{equation}}
\end{proposition}
\qed

From \eqref{226} one immediately deduces that for a subsequence, still denoted by $\eps$, one has 
\begin{equation}
\label{227}
\chi_{_{\oeps}}\rightarrow 1 \mbox{ a.e. in } \Omega \mbox{ as } \eps \rightarrow 0.
\end{equation}
\begin{proof}
In view of \eqref{57ter} one has 
\begin{equation}
\label{5002}
w^\eps \chi_{_{\oeps}}=w^\eps  \mbox{ a.e. in}\quad \Omega.
\end{equation}
Since $0\leq  \chi_{_{\oeps}}\leq 1$, one can extract a subsequence such that 
\begin{equation}
\label{5003}
\chi_{_{\oeps}}\rightharpoonup \theta \mbox{ in } L^\infty(\Omega) \mbox{ weakly-star as } \eps\rightarrow 0,
\end{equation}
so that using \eqref{cond3} and passing to the limit in \eqref{5002}, one has
$$
\theta=1,
$$
which implies that \eqref{5003} holds true with $\theta=1$ for the whole sequence $\eps$.

 It is then sufficient to write that
\begin{align*}
\|\chi_{_{\oeps}}-1\|_{L^1(\Omega)}=\into |\chi_{_{\oeps}}-1|=\into (1-\chi_{_{\oeps}})\rightarrow \into (1-\theta)=0 \mbox{ as }  \eps \to 0 
\end{align*}
 to deduce \eqref{226} for the whole sequence $\eps$.
 \end{proof}
 
 \medskip
 \begin{remark}
 \label{66661}
 Note that \eqref{227} implies that, for every subsequence $\eps'$ of $\eps$ and for almost every $x_0\in\Omega$, there exists $\eps_0(x_0)>0$ such that 
 $$
 \chi_{_{\Omega^{\eps'}}}(x_0)= 1 \,\, \forall \eps',\, \eps'<\eps_0(x_0),
 $$
 or in other terms that 
 \begin{equation}
 \label{520bis}
 x_0\in \Omega^{\eps'} \,\, \forall \eps',\, \eps'<\eps_0(x_0).
 \end{equation}
 \qed
\end{remark}
\medskip
\begin{remark}
\label{66662}
 Assertion \eqref{520bis} implies that almost every $x_0\in\Omega$ belongs to $\Omega^{\eps'}$ for $\eps'$ sufficiently small ($\eps'<\eps_0(x_0)$), which formally means that ``$\oeps$ is very close to $\Omega$". 
 
 In contrast, note that if we consider the case of holes periodically distributed at the vertices of a cubic lattice of size $\eps_j=1/2^j$,  with $j\in\mathbb{N}$, namely the case considered in the model example described in the Section~\ref{assumptions} above, every point of the form 
\begin{align*}
\dys c_k^{\eps_{j_0}}=k\,\eps_{j_0}=\left(\frac{k_1}{2^{j_0}}, \frac{k_2}{2^{j_0}},\cdots,\frac{k_N}{2^{j_0}}\right)
\dys\mbox{ with } k=(k_1,k_2,\cdots,k_N)\in\mathbb{Z}^N \mbox{ and } j_0\in\mathbb{N}
\end{align*}
is, for every $j\geq j_0$, the center of some hole $T_i^{\eps_j}$ which is extremely small, since its size is $r^{\eps_j}=C_0 \eps_j^{N/(N-2)}=C_0 2^{-jN/(N-2)}$; therefore such a point $c_k^{\eps_{j_0}}$ does not belong to $\oeps$ for $\eps_j=1/2^j\leq 1/2^{j_0}$; note that these points are dense in $\Omega$, and ``tend to invade the whole of $\Omega$" as $\eps_j$ tends to zero.
\qed
\end{remark}
\bigskip

\section{A priori estimates for the solutions\\ to the singular semilinear problem in $\oeps$}
\label{section7}
\bigskip
In this Section we state a priori estimates which are satisfied by every solution to \eqref{eqprimad} in the sense of Definition~\ref{sol}.

\begin{proposition}{\rm({\bf A priori estimate of $G_k(u^\eps)$ in $H_0^1(\oeps)$)}\break (Proposition~5.1 of \cite{GMM2})
\label{prop1}
Assume that the matrix $A$ and the function $F$ satisfy \eqref{eq0.0}, \eqref{car} and \eqref{eq0.1}. Then for every $u^\eps$ solution to problem \eqref{eqprimad} in the sense of Definition~\ref{sol} one has  
\begin{align}
\label{num3}
\|G_k(u^\eps)\|_{H_0^1(\oeps)}=\|DG_k(u^\eps)\|_{(L^2(\oeps))^N}\leq \frac{C_S}{\alpha} \frac{ \|h\|_{L^r(\oeps)}}{\Gamma(k)}\,\, \forall k>0,
\end{align}
where $C_S$ is the (generalized) Sobolev's constant defined by \eqref{2334}.}
\end{proposition}
\qed

\begin{remark} (Remark~5.2 of \cite{GMM2})
\label{rem52u}
From Poincar\'e's inequality 
\begin{equation}
\label{pis}
\|y^\eps\|_{L^2(\oeps)}\leq C_P(\oeps)\|Dy^\eps\|_{(L^2(\oeps))^N}\quad \forall y^\eps\in H_0^1(\oeps),
\end{equation}
where the constant $C_P(\oeps)$ is bounded independently of $\oeps$ when $\oeps\subset Q$, for $Q$ a bounded open set of $\mathbb{R}^N$, one deduces from \eqref{num3}, writing $u^\eps=T_k(u^\eps)+G_k(u^\eps)$, that every solution $u^\eps$ to problem \eqref{eqprimad} in the sense of Definition~\ref{sol} satisfies the following a priori estimate in $L^2(\oeps)$
\begin{align}
\label{53ter2}
\|u^\eps\|_{L^2(\oeps)}\dys\leq k|\oeps|^{\frac 12}+C_P(\oeps)\frac{C_S}{\alpha} \frac{ \|h\|_{L^r(\oeps)}}{\Gamma(k)}\,\, \forall k>0,
\end{align}
which, taking $k=k_0$ for some $k_0$ fixed or minimizing in $k$ provides an a priori estimate of $\|u^\eps\|_{L^2(\oeps)}$ which does not depend on $k$.
\qed
\end{remark}

\begin{proposition}{\rm({\bf A priori estimate of $\varphi^\eps DT_k(u^\eps)$ in $(L^2(\oeps))^N$ for\break $\varphi^\eps \in H_0^1(\oeps)\cap L^\infty(\oeps)$)} (Proposition~5.4 of \cite{GMM2})
\label{prop2}
Assume that the matrix $A$ and the function $F$ satisfy \eqref{eq0.0}, \eqref{car} and \eqref{eq0.1}. Then for every $u^\eps$ solution to problem \eqref{eqprimad} in the sense of Definition~\ref{sol} one has
\begin{align}
\label{619bis}
\begin{cases}
\dys\|\varphi^\eps DT_k(u^\eps)\|^2_{(L^2(\oeps))^N}\dys\leq \\\dys\leq\frac{32 k^2} {\alpha^2} \|A\|_{(L^\infty(\oeps))^{N\times N}}^2\|D\varphi^\eps\|_{(L^2(\oeps))^N}^2+\frac{C_S^2}{\alpha^2}\frac{\|h\|^2_{L^{r}(\oeps)}}{\Gamma(k)^2}\|\varphi^\eps\|_{L^\infty(\oeps)}^2\\
\dys \forall k>0,\,\, \forall\varphi^\eps\in H_0^1(\oeps)\cap L^\infty(\oeps),
\end{cases}
\end{align}
where $C_S$ is the (generalized) Sobolev's constant defined by \eqref{2334}.}
\end{proposition}
\qed
\begin{remark} (Remark~5.5 of \cite{GMM2})
\label{54bis}
From the a priori estimate \eqref{619bis} one deduces that every solution $u^\eps$ to problem \eqref{eqprimad} in the sense of Definition~\ref{sol} satisfies the following a priori estimate of  $\varphi^\eps T_k(u^\eps)$ in $H_0^1(\oeps)$
\begin{align}
\label{516bis}
\begin{cases}
\dys\|\varphi^\eps T_k(u^\eps)\|^2_{H_0^1(\oeps)}=\|D(\varphi^\eps T_k(u^\eps))\|^2_{(L^2(\oeps))^N}\leq\\
\dys \leq \left(\frac{64 k^2} {\alpha^2} \|A\|_{(L^\infty(\oeps))^{N\times N}}^2+2k^2\right) \|D\varphi^\eps\|^2_{(L^2(\oeps))^N}  + 2 \frac{C_S^2}{\alpha^2}\frac{\|h\|^2_{L^{r}(\oeps)}}{\Gamma(k)^2}\|\varphi^\eps\|_{L^\infty(\oeps)}^2\\
\dys \forall k>0,\,\, \forall\varphi^\eps\in H_0^1(\oeps)\cap L^\infty(\oeps).
\end{cases}
\end{align}
\qed
\end{remark}

For $\delta>0$,  define the function $Z_\delta: s\in[0,+\infty[\rightarrow Z_\delta(s)\in[0,+\infty[$ by 
\begin{equation}
\label{num23bis}
Z_\delta(s)=\begin{cases}
1& \mbox{if } 0\leq s\leq \delta, \\
 -\frac{s}{\delta}+2 & \mbox{if }  \delta\leq s\leq 2\delta ,\\
  0 & \mbox{if } 2\delta\leq s.
\end{cases}
\end{equation}

\begin{proposition}{\rm({\bf  Control of the quantity $\dys\into F(x,u^\eps)Z_\delta(u^\eps) v$ when $\delta$ is small)} (Proposition~5.9 of \cite{GMM2})
\label{prop3}
Assume that the matrix $A$ and the function $F$ satisfy \eqref{eq0.0}, \eqref{car} and \eqref{eq0.1}. Then for every $u^\eps$ solution to problem \eqref{eqprimad} in the sense of Definition~\ref{sol} and for every $v^\eps$ such  that
\begin{equation}
\label{5701}
\begin{cases}
\displaystyle v^\eps \in \mathcal V(\oeps),\,v^\eps\geq 0,\\ \dys\mbox{\rm{with }} -div \, {}^t\!A(x)Dv^\eps=\sum_{i \in I^\eps} \hat{ \varphi^\eps_i} (-div\, \hat{ g^\eps_i})+\hat{ f^\eps} \mbox{\rm{ in }} \mathcal{D}'(\oeps)\\
\mbox{\rm{where }} \hat{\varphi^\eps_i}\in H_0^1(\oeps)\cap L^\infty(\oeps), \hat{g^\eps_i}\in L^2(\oeps)^N, \hat{f^\eps}\in L^1(\oeps),\end{cases}
\end{equation}
one has
\begin{align}
\label{5701bis}
\begin{cases}
\dys\forall \delta>0,\,\,\int_{\oeps} F(x,u^\eps)\,Z_\delta(u^\eps)\, v^\eps\leq\\
\dys\leq \frac 32 \left(\int_{\oeps} \left|\sum_{i\in I^\eps} \hat{g^\eps_i} D\hat{\varphi^\eps_i}+ \hat{f^\eps}\right| \right)\delta+\sum_{i\in I^\eps} \into Z_\delta(u^\eps) \, \hat{g^\eps_i}Du^\eps \,\hat{\varphi^\eps_i}.
\end{cases}
\end{align}}
\end{proposition}
\qed

Note that the second term of the right-hand side of \eqref{5701bis} has a meaning since $Du^\eps\,\hat{\varphi}_i^\eps\in (L^2(\oeps))^N$ in view of \eqref{eq:phiDu}.

A consequence of Proposition~\ref{prop3} is:
\begin{proposition} {\rm (Proposition~5.12 of \cite{GMM2})
\label{prop69}
Assume that the matrix $A$ and the function $F$ satisfy \eqref{eq0.0}, \eqref{car} and \eqref{eq0.1}. Then for every $u^\eps$ solution to problem \eqref{eqprimad} in the sense of Definition~\ref{sol} one has
\begin{equation}
\label{5800}
\int_{\{u^\eps=0\}}F(x,u^\eps)v^\eps=0\quad \forall v^\eps\in\mathcal{V}(\oeps),v^\eps\geq 0.
\end{equation}}
\end{proposition}
\qed

%

 \smallskip

\section{Proof of the homogenization Theorem~\ref{homogenization}}
\label{section8}

\bigskip
\noindent{\bf First step}. In this step we state a priori estimates and we extract a subsequence still denoted by $\eps$ such that convergences \eqref{num2h}, \eqref{53bis} and \eqref{53ter} of Theorem~\ref{homogenization} hold true for some $u^0$ which satisfies \eqref{sol1h0}.

 As already said in the Existence Theorem~\ref{EUS} of Subsection~3.2, there exists at least one solution to problem \eqref{eqprimad} in the sense of Definition~\ref{sol}. This solution in particular satisfies the a priori estimates \eqref{num3}, \eqref{53ter2} and \eqref{516bis} stated in Proposition~\ref{prop1} and in Remarks~\ref{rem52u} and \ref{54bis} above.

Since $\oeps\subset\Omega$, since the generalized Sobolev's constant $C_S$ which appears in \eqref{2334} does not depend on $\oeps$ when $N\geq 3$, and is bounded independently of $\oeps$ when $N=2$ since $\oeps\subset \Omega$ (see the comment after \eqref{2334}), and since the Poincar\'e's constant $C_P(\oeps)$ which appears in \eqref{pis} is bounded independently of $\oeps$ since $\oeps\subset \Omega$ (see the comment after \eqref{pis}), the a priori estimates \eqref{num3} and \eqref{53ter2} imply that 
\begin{equation}
\label{APEG}
\|G_k(\reallywidetilde{u^\eps})\|_{H_0^1(\Omega)}= \|G_k(u^\eps)\|_{H_0^1(\oeps)}\leq C(k),
\end{equation}
\begin{equation}
\label{APEL}
\|\reallywidetilde{u^\eps}\|_{L^2(\Omega)}= \|u^\eps\|_{L^2(\oeps)}\leq C,
\end{equation}
where the constants $C(k)$ and $C$ do not depend on $\eps$ for $k>0$ fixed.

Similarly, taking in \eqref{516bis} $\varphi^\eps=z^\eps \varphi $, with $\varphi \in H_0^1(\Omega)\cap L^\infty(\Omega)$ and $z^\eps$ defined by Proposition~\ref{zeps}, and observing that $\|z^\eps \varphi\|_{L^\infty(\oeps)}$ and $\|D(z^\eps \varphi)\|_{(L^2(\Omega))^N}$ are bounded independently of $\eps$, one obtains that 
\begin{equation}
\label{APET2}
\|z^\eps \varphi T_k(\reallywidetilde{u^\eps})\|_{H_0^1(\Omega)}= \|z^\eps \varphi T_k(u^\eps)\|_{H_0^1(\oeps)}\leq C(k,\varphi),
 \end{equation}
 where the constant $C(k,\varphi)$ does not depend on $\eps$ for $k>0$ and $\varphi\in H_0^1(\Omega)\cap L^\infty(\Omega)$ fixed. 
 
 Using the (generalized) Sobolev's inequality \eqref{2334} for $G_k(\reallywidetilde{u^\eps})$, the fact that  $z^\eps$ is bounded in $H^1(\Omega)\cap L^\infty(\Omega)$, $\varphi\in H_0^1(\Omega)\cap L^\infty(\Omega)$ and \eqref{APEG}, one obtains that 
 \begin{equation}
 \label{73bis}
\|z^\eps \varphi G_k(\reallywidetilde{u^\eps})\|_{W_0^{1,q}(\Omega)}\leq C(k,\varphi)
 \dys\mbox{ where } q \mbox{ is defined by }  \frac{1}{q}=\frac{1}{2^*}+\frac{1}{2},
 \end{equation}
where $2^*$ is defined by \eqref{2333} and \eqref{212bis}.

Collecting together \eqref{APET2} and \eqref{73bis}  implies that $$ z^\eps\varphi \reallywidetilde{u^\eps}=z^\eps \varphi T_k( \reallywidetilde{u^\eps})+z^\eps \varphi G_k( \reallywidetilde{u^\eps})$$ is bounded in $W_0^{1,q}(\Omega)$, and therefore that
\begin{equation}
\label{74bis}
z^\eps\varphi \reallywidetilde{u^\eps} \mbox{ is compact in } L^q(\Omega) \mbox{ for every  } \varphi \in H_0^1(\Omega)\cap L^\infty(\Omega).
 \end{equation}
 
 On the other hand, let us write, for every $\varphi \in H_0^1(\Omega)\cap L^\infty(\Omega)$
\begin{equation}
\label{85bis}
 \varphi \reallywidetilde{u^\eps}=z^\eps\varphi   \reallywidetilde{u^\eps}+(1-z^\eps)\varphi  \reallywidetilde{u^\eps} .
\end{equation}
Since $(1-z^\eps)$ tends to zero in $L^p(\Omega)$ strongly for every $p<+\infty$ (see \eqref{10103} and \eqref{10102}), since $\varphi \in L^\infty(\Omega)$ and since $\reallywidetilde{u^\eps}$ is bounded in $L^2(\Omega)$ (see \eqref{APEL}), one has
 \begin{equation}
 \label{74ter}
 (1-z^\eps) \varphi  \reallywidetilde{u^\eps} \rightarrow 0 \mbox{ in } L^q(\Omega) \mbox{ strongly as } \eps\rightarrow 0.
 \end{equation}
  From \eqref{74bis}, \eqref{85bis} and \eqref{74ter} one concludes that 
 \begin{equation}
\label{APEC}
\varphi \reallywidetilde{u^\eps}\quad  \mbox{is compact in} \,\, L^q(\Omega) \quad \forall \varphi \in H_0^1(\Omega)\cap L^\infty(\Omega).
 \end{equation}
 
 In view of \eqref{APEL} and \eqref{APEC} one can extract a subsequence, still denoted by $\eps$, such that there exists some $u^0\in L^2(\Omega)$ such that  
 \begin{equation}
 \label{CV2}
 \reallywidetilde{u^\eps}\rightharpoonup u^0  \mbox{ in } L^2(\Omega)  \mbox{ weakly and a.e. in } \Omega \mbox{ as }  \eps\rightarrow 0.
 \end{equation}
 This proves \eqref{num2h}.
 
 For the same subsequence, one has, in view of \eqref{CV2}, \eqref{APEG}, \eqref{10102} and \eqref{APET2} 
  \begin{equation}
 \label{CV3}
G_k(\reallywidetilde{u^\eps}) \rightharpoonup G_k(u^0)  \mbox{ in }  H_0^1(\Omega)  \mbox{ weakly }\,\, \forall k>0 \mbox{ as }  \eps\rightarrow 0,
 \end{equation}
\begin{align}
 \label{CV4}
 \begin{cases}
z^\eps\varphi  T_k(\reallywidetilde{u^\eps}) \rightharpoonup \varphi T_k(u^0) \mbox{ in }  H_0^1(\Omega)  \mbox{ weakly}\\
 \forall k>0, \,\,   \forall \varphi \in H_0^1(\Omega)\cap L^\infty(\Omega), \mbox{ as }  \eps\rightarrow 0.
\end{cases}
 \end{align}
This proves \eqref{53bis}.

Moreover, since similarly to \eqref{APET2}, one has, taking $\varphi^\eps=w^\eps \varphi$ in \eqref{516bis}, 
\begin{equation}
\label{APET}
\|w^\eps \varphi  T_k(\reallywidetilde{u^\eps})\|_{H_0^1(\Omega)} \leq C(k,\varphi),
 \end{equation}
 where the constant $C(k,\varphi)$ does not depend on $\eps$ for  $k>0$ and $\varphi \in H_0^1(\Omega)\cap L^\infty(\Omega)$ fixed, one also has
 \begin{align}
 \label{CV5}
 \begin{cases}
w^\eps \varphi  T_k(\reallywidetilde{u^\eps}) \rightharpoonup \varphi T_k(u^0)  \mbox{ in }  H_0^1(\Omega)  \mbox{ weakly}\\
 \forall k>0, \,\,   \forall \varphi \in H_0^1(\Omega)\cap L^\infty(\Omega), \mbox{ as } \eps\rightarrow 0.
\end{cases}
 \end{align}
This proves \eqref{53ter}.

Note that since $\reallywidetilde{u^\eps}$ is nonnegative on $\Omega$, one has
\begin{equation}
\label{712bis}
u^0(x)\geq 0 \mbox{ a.e. } x\in\Omega.
\end{equation}
 
 On the other hand, since $G_k(u^0) \in H_0^1(\Omega)$ in view of \eqref{CV3} and since for every $\phi \in\mathcal{D}(\Omega)$ one has $\phi T_k(u^0)\in H_0^1(\Omega)$ in view of \eqref{CV4}, the function $u^0$ satisfies 
 \begin{equation}
\label{712ter}
 u^0 \in H^1_{\mbox{\tiny loc}}(\Omega).
\end{equation}
 
 As said in the introduction of this first step, we have extracted a subsequence and defined an $u^0$ which satisfies \eqref{sol1h0} such that convergences \eqref{num2h}, \eqref{53bis} and \eqref{53ter} hold true.
 
 \bigskip
\noindent{\bf Second step}.  We now consider any fixed $v\in\mathcal{V}(\Omega)$ with $v\geq 0$. In view of \eqref{10106} and of Remark~\ref{rem62}, the function $z^\eps v$ belongs to $\mathcal{V}(\oeps)$ with $z^\eps v\geq 0$ and satisfies \eqref{1021bis} when $v$ satisfies \eqref{10112}. The use of $v^\eps=z^\eps v$ in \eqref{condvh} is  therefore licit and one has 
\begin{align}
\label{1030}
\begin{cases}
\dys \int_{\oeps}\, {}^t\!A(x) D(z^\eps v)DG_k(u^\eps)\,+\\
\dys + \sum_{i\in I} \int_{\oeps} \hat{g_i} D(z^\eps \hat{\varphi_i}T_k(u^\eps))+\int_{\oeps} G^\eps D(w^\eps v T_k(u^\eps))\,+\\ \vspace{0.1cm}
\dys + \int_{\oeps}\left(z^\eps \hat{f}- \, {}^t\!A(x)Dv D z^\eps - \, {}^t\!A(x) Dz^\eps Dv\right)T_k(u^\eps)=\\ \vspace{0.1cm}
\dys =\langle-div\, {}^t\!A(x)D(z^\eps v), G_k(u^\eps)\rangle_{H^{-1}(\oeps),H_0^1(\oeps)}+ \langle\langle -div\, {}^t\!A(x)D(z^\eps v), T_k(u^\eps)\rangle\rangle_{\oeps}=\\
\dys= \int_{\oeps} F(x,u^\eps)z^\eps v.
\end{cases}
\end{align}
\indent From now on, $v\in\mathcal{V}(\Omega)$, $v\geq 0$, and $k>0$ will be fixed.
\smallskip

In the present step and in the next one, we pass to the limit, as $\eps$ tends to zero,  in the first term of the left-hand side of \eqref{1030} and we prove that
\begin{align}
\label{710bis}
\begin{cases}
\dys \int_{\oeps}\, {}^t\!A(x) D(z^\eps v)DG_k(u^\eps)\rightarrow\\
\dys \rightarrow \langle-div\, {}^t\!A(x)Dv, G_k(u^0)\rangle_{H^{-1}(\Omega),H_0^1(\Omega)} + \into G_k(u^0)\,v\,d\mu \mbox{ as } \eps\rightarrow 0.
\end{cases}
\end{align}
\indent For that we introduce, for $k>0$ fixed and for every $n>k$, the function\break 
$S_{k,n}:\mathbb{R}^+\rightarrow \mathbb{R}^+$ defined by
\begin{equation}
\label{defs}
S_{k,n}(s)=\begin{cases}
 0 & \mbox{if } 0\leq s\leq k, \\
 s-k & \mbox{if }  k\leq s\leq n,\\
n-k & \mbox{if } n\leq s.
\end{cases}
\end{equation}
Observe that one has
\begin{align}
\label{GSG}
G_k(s)=S_{k,n}(s)+G_n(s) \quad \forall s>0,\,\, \forall n,\,\, n>k,
\end{align}
\begin{align}
\label{STG}
S_{k,n}(s)=T_{n-k}(G_k(s)) \quad \forall s>0,\,\, \forall n,\,\, n>k.
\end{align}
\indent Using \eqref{GSG} we write the first  term of the left-hand side of \eqref{1030} as 
\begin{align}
\label{PalB}
\begin{cases}
\dys\int_{\oeps}\, {}^t\!A(x) D(z^\eps v)DG_k(u^\eps)=\\
\dys =\int_{\oeps}\, {}^t\!A(x) D(z^\eps v)DS_{k,n}(u^\eps)+\int_{\oeps}\, {}^t\!A(x) D(z^\eps v)DG_n(u^\eps).
\end{cases}
\end{align}
\indent We first pass to the limit in the first term of the right-hand side of \eqref{PalB} as $\eps$ tends to zero for $n$ and $k$ fixed, $n>k>0$. For that we write, using \eqref{10105} in the latest equality,
\begin{align}
\label{PalB1}
\begin{cases}
\dys\int_{\oeps}\, {}^t\!A(x) D(z^\eps v)DS_{k,n}(u^\eps)=\\\vspace{0.1cm}
\dys= \int_{\oeps}\, {}^t\!A(x) Dz^\eps D(S_{k,n}(u^\eps)v)- \int_{\oeps}\, {}^t\!A(x) Dz^\eps Dv\, S_{k,n}(u^\eps)\,+\\
\dys + \int_{\oeps}\, {}^t\!A(x) Dv\, DS_{k,n}(u^\eps)\,z^\eps =\\
\dys = \langle w^\eps \mu^\eps, S_{k,n}(u^\eps)v\rangle_{H^{-1}(\oeps),H_0^1(\oeps)}-\int_{\oeps}\, {}^t\!A(x) Dz^\eps Dv\, S_{k,n}(u^\eps)\,+\\
\dys + \int_{\oeps}\, {}^t\!A(x)  Dv\, DS_{k,n}(u^\eps)\, z^\eps.
\end{cases}
\end{align}

We now observe that in view of the convergence \eqref{CV3} of $G_k(\reallywidetilde{u^\eps})$ to $G_k(u^0)$ in $H_0^1(\Omega)$ weakly and of  formula \eqref{STG}, one has for $n>k$ fixed, 
\begin{align}
\label{719bis}
\begin{cases}
\dys S_{k,n}(\reallywidetilde{u^\eps})=T_{n-k}(G_k(\reallywidetilde{u^\eps}))\rightarrow T_{n-k}(G_k(u^0))=S_{k,n}(u^0)\\
\dys \mbox{in } H_0^1(\Omega) \mbox{ weakly and in } L^\infty(\Omega) \mbox{ weakly-star as } \eps\rightarrow 0.
\end{cases}
\end{align}

Therefore, using in the first term of the right-hand side of \eqref{PalB1} the strong convergence of $\mu^\eps$ to $\mu$ in $H^{-1}(\Omega)$ \big(see the fourth assertion of \eqref{cond5} and the convergence \eqref{cond3}\big), and then the equality \eqref{57bis2}, we have
\begin{align}
\label{PalB11}
\begin{cases}
\dys \langle w^\eps \mu^\eps, S_{k,n}(u^\eps)v \rangle_{H^{-1}(\oeps), H_0^1(\oeps)}
=\langle\mu^\eps, S_{k,n}(u^\eps)v  w^\eps  \rangle_{H^{-1}(\oeps), H_0^1(\oeps)}=\\
\dys = \langle\mu^\eps, S_{k,n}(\reallywidetilde{u^\eps})v  w^\eps  \rangle_{H^{-1}(\Omega), H_0^1(\Omega)}\rightarrow\\
\dys \rightarrow \langle\mu, S_{k,n}(u^0)v   \rangle_{H^{-1}(\Omega), H_0^1(\Omega)}= \into S_{k,n}(u^0)\,v\,d\mu\quad \mbox{as } \eps \rightarrow 0.
\end{cases}
\end{align}

For what concerns the second and the third terms of the right-hand side of \eqref{PalB1}, we have, in view of \eqref{10102} and \eqref{719bis}, 
\begin{align}
\label{PalB12}
\dys-\int_{\oeps}\, {}^t\!A(x) Dz^\eps Dv \,S_{k,n}(u^\eps)=-\int_{\Omega}\, {}^t\!A(x) Dz^\eps Dv \,S_{k,n}(\reallywidetilde{u^\eps})\rightarrow 0
\dys  \mbox{ as } \eps\rightarrow 0,
\end{align}
\begin{align}
\label{PalB13}
\begin{cases}
\dys\int_{\oeps}\, {}^t\!A(x) Dv DS_{k,n}(u^\eps)\,z^\eps=\int_{\Omega}\, {}^t\!A(x) Dv DS_{k,n}(\reallywidetilde{u^\eps})\,z^\eps \rightarrow \\
\dys \rightarrow \int_{\Omega}\, {}^t\!A(x) Dv DS_{k,n}(u^0) \quad \mbox{ as } \eps\rightarrow 0.
\end{cases}
\end{align}

Collecting together \eqref{PalB1}, \eqref{PalB11}, \eqref{PalB12} and \eqref{PalB13}, we have proved that the first term of the right-hand side of \eqref{PalB} satisfies 
\begin{align}
\label{PalB1n}
\begin{cases}
\dys\int_{\oeps}\, {}^t\!A(x) D(z^\eps v)DS_{k,n}(u^\eps)\rightarrow\\
\dys \rightarrow \int_{\Omega}\, {}^t\!A(x) DvDS_{k,n}(u^0)+\int_{\Omega}\, S_{k,n}(u^0)\,v\,d\mu \quad \mbox{ as } \eps\rightarrow 0.
\end{cases}
\end{align}

Let us now pass to the limit in the right-hand side of \eqref{PalB1n} as $n$ tends to infinity.

Since $DS_{k,n}(u^0)=DG_k(u^0)\chi_{_{\{k\leq u^0\leq n\}}}$, one has
$$
S_{k,n}(u^0)\rightarrow G_k(u^0) \quad \mbox{in } H_0^1(\Omega) \mbox{ strongly as } n\rightarrow +\infty,
$$
and therefore
$$
S_{k,n}(u^0)\rightarrow G_k(u^0) \quad \mbox{in } L^1(\Omega;d\mu) \mbox{ strongly as } n\rightarrow +\infty.
$$

Therefore the right-hand side of \eqref{PalB1n} satisfies, as $n$ tends to infinity, since $v\in L^\infty(\Omega;d\mu)$ (see \eqref{58bis}),
\begin{align}
\label{724bis}
\begin{cases}
\dys \int_{\Omega}\, {}^t\!A(x) DvDS_{k,n}(u^0)+\int_{\Omega}\, S_{k,n}(u^0)\,v\,d\mu\rightarrow\\
\dys  \rightarrow\int_{\Omega}\, {}^t\!A(x) DvDG_k(u^0)+\int_{\Omega}\, G_k(u^0)\,v\,d\mu=\\
\dys = \langle  -div\, {}^t\!A(x)Dv,G_k(u^0)\rangle_{H^{-1}(\Omega),H_0^1(\Omega)}+\into G_k(u^0)\,v\,d\mu \mbox{ as } n\rightarrow +\infty.
\end{cases}
\end{align}

Passing to the limit in \eqref{PalB} first as $\eps$ tends to zero for $n$ fixed and then for $n$ tending to infinity, and collecting together \eqref{PalB1n} and \eqref{724bis} will prove \eqref{710bis} whenever we will have proved that the second term of the right-hand side of \eqref{PalB} satisfies 
\begin{align}
\label{estimfinal2}
\dys \limsup_{\eps} \left| \int_{\oeps}\, {}^t\!A(x) D(z^\eps v)DG_n(u^\eps)\right|\rightarrow 0
\dys \mbox{ as } n\to +\infty,
\end{align}
see the Third step just below.

\bigskip
\noindent{\bf Third step}. In this step we prove \eqref{estimfinal2}. As just said, this will complete the proof of \eqref{710bis}. For that, we estimate the second term of the right-hand side of \eqref{PalB}.

 Since $z^\eps$ is bounded in $H^1(\Omega)\cap L^\infty(\Omega)$ (see \eqref{10103} and \eqref{10102}), one has
\begin{align}
\label{estim1}
\begin{cases}
\dys\left|\int_{\oeps}\, {}^t\!A(x) D(z^\eps v)DG_n(u^\eps)\right|\leq\\
\dys\leq \|A\|_{(L^\infty(\Omega))^{N\times N}}\|D(z^\eps v)\|_{(L^2(\Omega))^N}\|DG_n(u^\eps)\|_{(L^2(\oeps))^N}\leq\\
\dys \leq  \|A\|_{(L^\infty(\Omega))^{N\times N}}\left(\|z^\eps\|_{L^\infty(\Omega)}\|Dv\|_{(L^2(\Omega))^N}+\|v\|_{L^\infty(\Omega)}\|Dz^\eps\|_{_{(L^2(\Omega))^N}}\right)\\
\hspace{8cm}\dys\|DG_n(u^\eps)\|_{(L^2(\oeps))^N}\leq\\
\dys \leq C(v)\|DG_n(u^\eps)\|_{(L^2(\oeps))^N} \quad \forall \eps, \quad \forall n,
\end{cases}
\end{align}
where $C(v)$ is a constant which depends on $v$ but neither on $\eps$ nor on $n$. 

We now estimate $\|DG_n(u^\eps)\|_{(L^2(\oeps))^N}$ in a way which is more precise than the a priori estimate \eqref{num3}. For that we use the (energy) equality (5.4) of \cite{GMM2}, namely 
\begin{equation*}
                          \int_{\oeps} A(x) DG_n(u^\eps) D G_n(u^\eps)=\int_{\oeps}  F(x,u^\eps)G_n(u^\eps),
\end{equation*}
which is formally obtained by using $G_n(u^\eps)$ as test function in \eqref{eqprimad}. Using in this inequality the coercivity \eqref{eq0.0} of the matrix $A$ and the growth condition \eqref{eq0.1} on  the function $F$ gives, since $\Gamma$ is increasing and since $G_n(s)=0$ for $s\leq n$,
\begin{align*}
\dys\alpha \|DG_n(u^\eps)\|^2_{(L^2(\oeps))^N}\leq \int_{\oeps} \frac{h(x)}{\Gamma(u^\eps)}G_n(u^\eps)\leq \int_{\oeps} \frac{h(x)}{\Gamma(n)}G_n(u^\eps)= \into \frac{h(x)}{\Gamma(n)}G_n(\reallywidetilde{u^\eps}).
\end{align*}

Passing to the limit in $\eps$  for $n$ fixed thanks to \eqref{CV3} gives
\begin{equation}
\label{estimG_k}
\limsup_{\eps} \|DG_n(u^\eps)\|^2_{(L^2(\oeps))^N}\leq \omega(n) \quad \forall n>0,
\end{equation}
where $\omega(n)$ is defined by
\begin{equation}
\label{727bis}
\omega^2(n)=\frac{1}{\alpha}\into \frac{h(x)}{\Gamma(n)}G_n(u^0)\quad \forall n>0.
\end{equation}

Since $\Gamma$ is increasing and since for $s\geq 0$ fixed $G_n(s)$ is nonincreasing in $n$, one has, for $n\geq n_0$, 
\begin{align*}
\dys\omega^2(n)= \frac{1}{\alpha}\into \frac{h(x)}{\Gamma(n)}G_n(u^0)\chi_{_{\{u^0\geq n\}}}
\dys\leq \frac{1}{\alpha}\into \frac{h(x)}{\Gamma(n_0)}G_{n_0}(u^0)\chi_{_{\{u^0\geq n\}}} \quad \forall n,\, n\geq n_0.
\end{align*}

Since the measure of the set $\{x\in\Omega: u^0(x)\geq n\}$ tends to zero as $n$ tends to infinity (recall that $u^0\in L^2(\Omega)$), and since $h(x)G_{n_0}(u^0)\in L^1(\Omega)$, one deduces, fixing $n_0$, that
\begin{align}
\label{727ter}
\dys\omega^2(n)\rightarrow 0
\dys\mbox{ as } n\to +\infty.
\end{align}

Collecting together \eqref{estim1}, \eqref{estimG_k}, \eqref{727bis} and \eqref{727ter} proves that the second term of the right-hand side of \eqref{PalB} satisfies 
\begin{align*}
\dys \limsup_{\eps} \left| \int_{\oeps}\, {}^t\!A(x) D(z^\eps v)DG_n(u^\eps)\right|\leq C(v)\omega(n)\rightarrow 0,
\dys \mbox{ as } n\to +\infty.
\end{align*}
i.e. \eqref{estimfinal2}.

\bigskip
\noindent{\bf Fourth step}. In this step we pass to the limit, as $\eps$ tends to zero, in the second, third and fourth terms of the left-hand side of \eqref{1030} and we prove that 
\begin{align}
\label{7000}
\begin{cases}
\dys \sum_{i\in I} \int_{\oeps} \hat{g_i} D(z^\eps \hat{\varphi_i}T_k(u^\eps))+\int_{\oeps} G^\eps D(w^\eps v T_k(u^\eps))\,+\\
\dys + \int_{\oeps}\left(z^\eps \hat{f}- \, {}^t\!A(x)Dv D z^\eps - \, {}^t\!A(x) Dz^\eps Dv\right)\,T_k(u^\eps)\rightarrow\\
\dys \rightarrow \langle \langle-div\, {}^t\!A(x)Dv, T_k(u^0)\rangle\rangle_\Omega+\into T_k(u^0)\,v\,d\mu\quad \mbox{as } \eps\to0,
\end{cases}
\end{align}
where in the last line we used the notation \eqref{dc} of Definition~\ref{def32}.

For the second term of the left-hand side of \eqref{1030}, we have, in view of \eqref{CV4} and since $\hat{\varphi}_i\in H_0^1(\Omega)\cap L^\infty(\Omega)$, 
\begin{align}
\label{1032}
\dys \int_{\oeps} \hat{g_i} D(z^\eps \hat{\varphi_i}T_k(u^\eps))= \int_{\Omega} \hat{g_i}  D(z^\eps \hat{\varphi_i}T_k(\reallywidetilde{u^\eps}))
\dys \rightarrow  \into \hat{g_i}  D(\hat{\varphi_i}T_k(u^0))\quad \mbox{as } \eps\rightarrow 0.
\end{align}

Similarly, using the strong convergence of $G^\eps$ to $G$ in $(L^2(\Omega))^N$ (see \eqref{10112bis}), \eqref{CV5} and the fact that $v\in H_0^1(\Omega)\cap L^\infty(\Omega)$, we have, for the third term of the left-hand side of \eqref{1030},
\begin{align}
\label{1033}
\dys\int_{\oeps} G^\eps D(w^\eps v T_k(u^\eps))=\into G^\eps D(w^\eps v T_k(\reallywidetilde{u^\eps})) 
\dys\rightarrow \into G D(v T_k(u^0)) \quad \mbox{as } \eps\rightarrow 0.
\end{align}
Moreover, in view of \eqref{10112bis} and of \eqref{57bis2} we have, since $vT_k(u^0)\in H_0^1(\Omega)\cap L^\infty(\Omega)$ because of \eqref{CV4},
\begin{align}
\label{1033bis}
\begin{cases}
\dys\into G D(vT_k(u^0))=\langle-div\, G,vT_k(u^0)\rangle_{H^{-1}(\Omega),H_0^1(\Omega)}=\\
\dys =\langle\mu,vT_k(u^0)\rangle_{H^{-1}(\Omega),H_0^1(\Omega)}=\into T_k(u^0) \,v\, d\mu.
\end{cases}
\end{align}

Finally, in view of \eqref{10102} and \eqref{CV2}, and using the fact that $T_k(\reallywidetilde{u^\eps})Dv$ converges to $T_k(u^0)Dv$ in $(L^2(\Omega))^N$ strongly by Lebesgue's dominated convergence theorem, we have, for the fourth term of the left-hand side of \eqref{1030},
\begin{align}
\label{1034}
\begin{cases}
\dys\int_{\oeps} \left(z^\eps \hat{f}-\, {}^t\!A(x) Dv Dz^\eps-\, {}^t\!A(x) Dz^\eps Dv\right)\,T_k(u^\eps)=\\\vspace{0.1cm}
\dys = \into \left(z^\eps\hat{f}-\, {}^t\!A(x) Dv Dz^\eps-\, {}^t\!A(x) Dz^\eps Dv\right)\,T_k(\reallywidetilde{u^\eps}) \rightarrow \into \hat{f} T_k(u^0) \quad \mbox{as } \eps\rightarrow 0.
\end{cases}
\end{align}

Collecting together \eqref{1032}, \eqref{1033}, \eqref{1033bis} and \eqref{1034} we have proved that the second, third and fourth terms of the left-hand side of \eqref{1030} satisfy 
\begin{align}
\label{1034bis}
\begin{cases}
\dys \sum_{i\in I} \int_{\oeps} \hat{g_i} D(z^\eps \hat{\varphi_i}T_k(u^\eps))+\int_{\oeps} G^\eps D(w^\eps v T_k(u^\eps))\,+\\
\dys + \int_{\oeps}\left(z^\eps \hat{f}- \, {}^t\!A(x)Dv D z^\eps - \, {}^t\!A(x) Dz^\eps Dv\right)\,T_k(u^\eps)\rightarrow\\
\dys \rightarrow \sum_{i\in I} \into \hat{g_i} D( \hat{\varphi_i}T_k(u^0))+\into T_k(u^0)\,v\,d\mu+ \into \hat{f}T_k(u^0) \mbox{ as } \eps\to0.
\end{cases}
\end{align}
But in view of the notation \eqref{dc} of Definition~\ref{def32}, one has 
\begin{equation}
\label{729bis}
\dys \sum_{i\in I} \int_{\oeps} \hat{g_i} D(\hat{\varphi_i}T_k(u^0))+  \into \hat{f}T_k(u^0)= \langle \langle-div\, {}^t\!A(x)Dv,T_k(u^0)\rangle\rangle_\Omega.
\end{equation}
 
 From \eqref{1034bis} and \eqref{729bis} one deduces \eqref{7000}.
 \bigskip
 
\noindent{\bf Fifth step}.  At this point, see \eqref{710bis} and \eqref{7000}, we passed to the limit in the left-hand side of \eqref{1030}. In the sixth, seventh and eighth steps, we will pass to the limit in the right-hand side of \eqref{1030}. 

Before of that, we prove in the present step that 
\begin{equation}
\label{10400}
\into F(x,u^0)v<+\infty\quad \forall v\in\mathcal{V}(\Omega), v\geq 0,
\end{equation}
or in other terms that assertion (\ref{condv3h} {\it i}) holds true.

Since the left-hand side of \eqref{1030} converges as $\eps$ tends to zero, the right-hand side of \eqref{1030} satisfies
\begin{equation}
\label{437bis}
\int_{\oeps} F(x,u^\eps)z^\eps v\leq C(v)\quad \forall v\in\mathcal{V}(\Omega), v\geq 0, \,\, \forall \eps,
\end{equation}
the constant $C(v)<+\infty$ does not depend on $\eps$. Using the extension by zero defined in \eqref{defex}, \eqref{437bis} is equivalent to
\begin{equation}
\label{10402}
\into \reallywidetilde{F(x,u^\eps)}z^\eps v\leq C(v) \,\, \forall \eps.
\end{equation}

We claim that 
\begin{equation}
\label{10403}
\reallywidetilde{F(x,u^\eps)}\rightarrow F(x,u^0)\quad \mbox{a.e. } x\in\Omega \quad \mbox{as }  \eps\rightarrow 0.
\end{equation}
Indeed, in view of \eqref{520bis}, we know that, for every subsequence $\eps'$ of $\eps$ and for almost every $x_0\in\Omega$, there exists $\eps_0(x_0)$ such that $x_0$ belongs to $\Omega^{\eps'}$ for every $\eps'<\eps_0(x_0)$. This implies that 
$$
\reallywidetilde{F(x_0,u^{\eps'}(x_0))}=F(x_0,\reallywidetilde{u^{\eps'}}(x_0))\quad \forall \eps',\, \eps'<\eps_0(x_0).
$$
Since 
$$
F(x_0,\reallywidetilde{u^{\eps'}}(x_0))\rightarrow F(x_0, u^0(x_0)) \quad \mbox{a.e. } x_0\in\Omega
$$
in view of the convergence \eqref{CV2} and of the Carath\'eodory hypothesis \eqref{car}, this implies \eqref{10403}.

Results \eqref{10402} and \eqref{10403} combined with \eqref{10102}, the fact that $F(x,\reallywidetilde{u^\eps})z^\eps v\geq 0$, and finally Fatou's Lemma immediately imply \eqref{10400}.

\bigskip
\noindent{\bf Sixth step}. From now on, we introduce a new parameter $\delta>0$ and we write the right-hand side of \eqref{1030} as
\begin{align}
\label{10410}
\dys\int_{\oeps} F(x,u^\eps) z^\eps v= \int_{\oeps} F(x,u^\eps) \, Z_\delta(u^\eps)\,z^\eps v+\int_{\oeps} F(x,u^\eps)\, (1-Z_\delta(u^\eps))\, z^\eps v,
\end{align}
where $Z_\delta$ is the function defined by \eqref{num23bis}.

\bigskip
In the present step we prove that the first term of the right-hand side of \eqref{10410} satisfies
\begin{equation}
\label{1035}
\limsup_{\eps}\int_{\oeps}F(x,u^\eps)\, Z_\delta(u^\eps)\, z^\eps v\rightarrow0\quad  \mbox{as } \delta\to 0.
\end{equation}

For that we use estimate \eqref{5701bis} of Proposition~\ref{prop3} above with $v^\eps=z^\eps v$ for any $v\in\mathcal{V}(\Omega)$, $v\geq 0$; this choice is licit in view of \eqref{10106}. In view of \eqref{1021bis}, the estimate reads as

\begin{align}
\label{1040}
\begin{cases}
\dys\forall \delta>0,\, \int_{\oeps}F(x,u^\eps)\, Z_\delta(u^\eps)\, z^\eps v\leq I_\delta^\eps+ I\!I_\delta^\eps, \\
\dys \mbox{where}\\
\dys I_\delta^\eps= \frac32  \left(\int_{\oeps}\left|\sum_{i\in I}\hat{g_i}D(z^\eps \hat{\varphi_i})+G^\eps D(w^\eps v)\,+\right.\right.\\
\hspace{2cm}\dys  \left.\left.+z^\eps \hat{f}-\, {}^t\!A(x) Dv Dz^\eps-\, {}^t\!A(x) Dz^\eps Dv\right| \right)\delta,\\
\dys I\!I_\delta^\eps=  \sum_{i\in I} \int_{\oeps} Z_\delta(u^\eps)\, \hat{g_i}\,D u^\eps\, z^\eps \hat{\varphi_i}+ \int_{\oeps} Z_\delta(u^\eps) \,G^\eps Du^\eps\, w^\eps v.
\end{cases}
\end{align}

Since $z^\eps$ and $w^\eps$ are bounded in $H^1(\Omega)\cap L^\infty(\Omega)$ (see \eqref{10102} and \eqref{cond3}), and since $\hat{\varphi}_i$ and $v$  belong to $H_0^1(\Omega)\cap L^\infty(\Omega)$, since $\hat{f}\in L^1(\Omega)$ (see \eqref{10112}) and since $G^\eps$ is bounded in $(L^2(\Omega))^N$ (see \eqref{10112bis}), we have, as far as $I_\delta^\eps$ is concerned, 
\begin{equation}
\label{1036bis}
\limsup_{\eps} I_{\delta}^\eps\leq \frac 32 C\delta,
\end{equation}
where the constant $C$ does not depend neither on $\eps$ nor on $\delta$,
and therefore we have
\begin{equation}
\label{1036ter}
\limsup_{\eps} I_{\delta}^\eps\rightarrow 0 \quad \mbox{as } \delta \to 0.
\end{equation}

For what concerns $I\!I_\delta^\eps$, we write for the first term
\begin{align*}
\begin{cases}
\dys D u^\eps \, z^\eps \hat{\varphi_i}=(DT_k(u^\eps)+DG_k(u^\eps))z^\eps \hat{\varphi_i}=\\
\dys= D (z^\eps \hat{\varphi_i}T_k(u^\eps))-T_k(u^\eps) D(z^\eps \hat{\varphi_i})+DG_k(u^\eps)\,z^\eps \hat{\varphi_i}\,\, \mbox{ in } \mathcal{D'}(\oeps).
\end{cases}
\end{align*}
Since $z^\eps \hat{\varphi_i} T_k(u^\eps)$ and $G_k(u^\eps)$ belong to $H_0^1(\oeps)$ (see \eqref{sol1h} and \eqref{10106}), we have in view of \eqref{57bis}
\begin{align}
\begin{cases}
\label{1036quarto}
\dys \reallywidetilde{D u^\eps\, z^\eps \hat{\varphi_i}}=D\reallywidetilde{u^\eps} z^\eps\varphi_i=
\\\dys= D (z^\eps\hat{\varphi_i}T_k(\reallywidetilde{u^\eps}))-T_k(\reallywidetilde{u^\eps}) D(z^\eps \hat{\varphi_i})+DG_k(\reallywidetilde{u^\eps})z^\eps\hat{\varphi_i} \,\, \mbox{ in } \mathcal{D}'(\Omega).
\end{cases}
\end{align}
Since each term of the right-hand side of \eqref{1036quarto} is zero on $\Omega\setminus\oeps$, and even if
$$
\reallywidetilde{Z_\delta(u^\eps)}=0\,\, \mbox{ in }\,\, \Omega\setminus \oeps \quad \mbox{while}\quad  Z_\delta(\reallywidetilde{u^\eps})=1 \mbox{ in } \Omega\setminus\oeps,
$$
we have
\begin{align}
\label{1036quinto}
\begin{cases}
\dys\reallywidetilde{Z_\delta(u^\eps)\hat{g_i}Du^\eps(z^\eps\varphi_i)}=Z_\delta(\reallywidetilde{u^\eps})\hat{g}_iD\reallywidetilde{u^\eps}\, z^\eps \hat{\varphi}_i=\\
\dys=Z_\delta(\reallywidetilde{u^\eps})\hat{g_i}	\left(D(z^\eps \hat{\varphi_i}T_k(\reallywidetilde{u^\eps}))-T_k(\reallywidetilde{u^\eps})D(z^\eps \hat{\varphi_i})+DG_k(\reallywidetilde{u^\eps})z^\eps \hat{\varphi_i}\right) \mbox{ in } \mathcal{D}'(\Omega).
\end{cases}
\end{align}
Therefore, in view of \eqref{CV2}, \eqref{CV3},  \eqref{CV4} and \eqref{10102}, we have
\begin{align}
\label{1042}
\begin{cases}
\dys \forall\delta>0 \mbox{ fixed},\\
\dys\int_{\oeps}Z_\delta(u^\eps)\, \hat{g_i}\,Du^\eps \,z^\eps \hat{\varphi_i}=\\\vspace{0.1cm}
\dys =\into Z_\delta(\reallywidetilde{u^\eps})\,\hat{g_i}	\left( D(z^\eps \hat{\varphi_i}T_k(\reallywidetilde{u^\eps}))-T_k(\reallywidetilde{u^\eps})D(z^\eps \hat{\varphi_i})+DG_k(\reallywidetilde{u^\eps})z^\eps \hat{\varphi_i}\right)\rightarrow\\\vspace{0.1cm}
\dys \rightarrow \into Z_\delta(u^0)\,\hat{g_i}\left(D(\hat{\varphi_i}T_k(u^0))-T_k(u^0)D\hat{\varphi_i}+DG(u^0)\hat{\varphi_i}\right)=\\
\dys= \into Z_\delta(u^0)\, \hat{g_i} \,Du^0\, \hat{\varphi_i}\quad \mbox{as } \eps\to 0.
\end{cases}
\end{align}

A proof which is very similar to the proof of \eqref{1042} implies that for the second term of $I\!I_\delta^\eps$ we have
\begin{align}
\label{1043}
\begin{cases}
\dys\forall\delta>0 \mbox{ fixed},\\
\dys \int_{\oeps} Z_\delta(u^\eps)\,G^\eps Du^\eps\, w^\eps v\rightarrow \into Z_\delta(u^0)\,G^0 Du^0 \,v\quad\mbox{ as } \eps\to 0.
\end{cases}
\end{align}

Let us now pass to the limit in the right-hand side of \eqref{1042} as $\delta$ tends to zero. Since
$$
Z_\delta(s)\rightarrow \chi_{_{\{s=0\}}}(s),\quad \forall s\geq 0,
$$
and since $u^0 \in H^1_{\mbox{\tiny loc}}(\Omega)$ implies that 
$$
Du^0=0\mbox{ a.e. in } \{x\in\Omega:u^0(x)=0\},
$$
we have
\begin{equation}
\label{1044}
\dys\into Z_\delta(u^0)\,\hat{g_i}\,Du^0\hat{\varphi_i}\rightarrow \into  \chi_{_{\{u^0=0\}}}\hat{g_i}\,Du^0\,\hat{\varphi_i}=0\quad \mbox{as } \delta\to 0.
\end{equation}

The same proof implies that for the right-hand side of \eqref{1043} we have
\begin{equation}
\label{1045}
\dys\into Z_\delta(u^0)\,G^0Du^0v\rightarrow 0\quad \mbox{as } \delta\to 0.
\end{equation}

Collecting together the definition \eqref{1040} of $I\!I^\eps_\delta$ and the results obtained in \eqref{1042},  \eqref{1044}, \eqref{1043} and \eqref{1045} proves that 
\begin{equation}
\label{1046}
\lim_{\eps} I\!I^\eps_\delta\rightarrow0 \quad\mbox{as  } \delta\to 0.
\end{equation}

Finally, collecting together \eqref{1040}, \eqref{1036ter} and \eqref{1046} proves \eqref{1035}.

\bigskip
\noindent{\bf Seventh step}. In this step we prove that 
\begin{equation}
\label{1050}
\int_{\Omega\cap\{u^0=0\}} F(x,u^0)\,v=0.
\end{equation}

Since by the definition \eqref{defex}, one has $\reallywidetilde{F(x,u^\eps)}=0$ on $\Omega\setminus\oeps$, and even if 
$$
\reallywidetilde{Z_\delta(u^\eps)}=0\,\, \mbox{ in }\,\, \Omega\setminus\oeps \quad \mbox{while} \quad Z_\delta(\reallywidetilde{u^\eps})=1 \mbox{ in } \Omega\setminus\oeps,
$$
 we have
\begin{equation}
\label{1043bis}
\reallywidetilde{F(x,u^\eps)Z_\delta(u^\eps)}=\reallywidetilde{F(x,u^\eps)} Z_\delta(\reallywidetilde{u^\eps}) \mbox{ in } \Omega,
\end{equation}
from which we deduce that 
\begin{align}
\label{75454}
\begin{cases}\vspace{0.1cm}
\dys \forall \delta>0, \,\, \int_{\Omega\cap \{u^0=0\}} \reallywidetilde{F(x,u^\eps)}\,Z_\delta(\reallywidetilde{u^\eps})\,z^\eps v=\int_{\oeps\cap \{u^0=0\}} F(x,u^\eps)\,Z_\delta(u^\eps)\,z^\eps v\leq\\
\dys \leq \int_{\oeps}F(x,u^\eps)\,Z_\delta(u^\eps)\,z^\eps v. 
\end{cases}
\end{align}

Since in view of \eqref{CV2} we have
$$
Z_\delta(\reallywidetilde{u^\eps})\rightarrow Z_\delta(u^0)=1 \mbox{ a.e. in } \{x\in\Omega: u^0(x)=0\} \mbox{ as } \eps\rightarrow 0,
$$
we have, in view of \eqref{10403} and \eqref{10102}
\begin{equation}
\label{75555}
\reallywidetilde{F(x,u^\eps)}Z_\delta(\reallywidetilde{u^\eps})z^\eps v\rightarrow F(x,u^0)\,v \mbox{ a.e. in } \{x\in\Omega: u^0(x)=0\} \mbox{ as } \eps\rightarrow 0.
\end{equation}

Using Fatou's Lemma in the left-hand side of \eqref{75555} we obtain 
\begin{align}
\label{75656}
\begin{cases}\vspace{0.1cm}
\dys \forall \delta>0, \,\, \int_{\Omega\cap \{u^0=0\}} F(x,u^0) v\leq \liminf_{\eps} \int_{\Omega\cap \{u^0=0\}} \reallywidetilde{F(x,u^\eps)}\,Z_\delta(\reallywidetilde{u^\eps})\,z^\eps v\leq\\
\dys \leq \limsup_{\eps} \int_{\oeps}F(x,u^\eps)\,Z_\delta(u^\eps)\,z^\eps v,
\end{cases}
\end{align}
which letting $\delta$ tend to zero and using \eqref{1035} implies \eqref{1050}.

\bigskip
\noindent{\bf Eight step}. In this step we prove that the second term of the right-hand side of \eqref{10410} satisfies 
\begin{align}
\label{75757}
\dys \lim_{\eps} \int_{\oeps} F(x,u^\eps)(1-Z_\delta(u^\eps))\,z^\eps v\rightarrow
\dys  \int_{\Omega}F(x,u^0)\, v \mbox{ as } \delta\rightarrow 0.
\end{align}
Indeed, similarly to the results obtained in the seventh step, we have
\begin{equation}
\label{75858}
\reallywidetilde{F(x,u^\eps)(1-Z_\delta(u^\eps))}=\reallywidetilde{F(x,u^\eps)}(1-Z_\delta(\reallywidetilde{u^\eps}))\mbox{ in } \Omega,
\end{equation}
as well as
\begin{align}
\label{75959}
\dys \reallywidetilde{F(x,u^\eps)}(1-Z_\delta(\reallywidetilde{u^\eps}))z^\eps v\rightarrow F(x,u^0)(1-Z_\delta(u^0))v \mbox{ a.e. in } \Omega.
\end{align}

On the other hand, we have 
$$
1-Z_\delta(\reallywidetilde{u^\eps})=0 \mbox{ a.e. in }  \{x\in \Omega:\reallywidetilde{u^\eps}(x)\leq \delta\},
$$
while in view of the conditions (\ref{eq0.1} {\it iii}) and (\ref{eq0.1} {\it ii}) on the functions $F(x,s)$ and $\Gamma(s)$, we have
$$
0\leq \reallywidetilde{F(x,u^\eps)}\leq \frac{h(x)}{\Gamma(\reallywidetilde{u^\eps})}\leq \frac{h(x)}{\Gamma(\delta)} \mbox{ a.e. in }  \{x\in \Omega:\reallywidetilde{u^\eps}(x)> \delta\}.
$$
Together with $0\leq Z_\delta(s)\leq 1$ and \eqref{10103}, this implies that 
\begin{equation}
\label{76060}
0\leq \reallywidetilde{F(x,u^\eps)} (1-Z_\delta(\reallywidetilde{u^\eps}))z^\eps v\leq \frac{h(x)}{\Gamma(\delta)}v \mbox{ a.e. in } \Omega,
\end{equation}
where $h(x)v\in L^1(\Omega)$ in view of condition (\ref{eq0.1} {\it i}).

From \eqref{75858}, \eqref{75959}, \eqref{76060} and Lebesgue's dominated convergence theorem we deduce that 
\begin{align}
\label{76161}
\begin{cases}
\dys \forall \delta>0, \,\, \int_{\oeps} F(x,u^\eps)(1-Z_\delta(u^\eps))\,z^\eps v=\int_{\oeps} \reallywidetilde{F(x,u^\eps)}\,(1-Z_\delta(\reallywidetilde{u^\eps}))\,z^\eps v\rightarrow\\\vspace{0.1cm}
\dys \rightarrow \into F(x,u^0)\,(1-Z_\delta(u^0)) \,v\,\mbox{ as } \eps\rightarrow 0.
\end{cases}
\end{align}

Since 
$$
Z_\delta(u^0)\rightarrow \chi_{_{\{u^0=0\}}}\mbox{ a.e. in } \Omega \mbox{ as } \delta \rightarrow 0,
$$
applying again Lebesgue's dominated convergence theorem and \eqref{1050} implies that 
\begin{align*}
\dys \int_{\Omega} F(x,u^0)(1-Z_\delta(u^0)) \,v\rightarrow \into F(x,u^0)(1-\chi_{_{\{u^0=0\}}})\,v=\into F(x,u^0)v \,\mbox{ as } \delta\rightarrow 0.
\end{align*}

This proves \eqref{75757}.

\bigskip
\noindent{\bf Ninth (and last) step}.  Collecting together \eqref{1030}, \eqref{710bis}, \eqref{7000}, \eqref{10410}, \eqref{1035} and \eqref{75757},  we have proved that $u^0$ satisfies (\ref{condv3h} {\it ii}). We also have proved in \eqref{10400} that $u^0$ satisfies (\ref{condv3h} {\it i}).

Since we have proved in the first step that the subsequence that we have extracted satisfies the convergences  \eqref{num2h}, \eqref{53bis} and \eqref{53ter}, and that $u^0$ satisfies \eqref{sol1h0}, the proof of Theorem~\ref{homogenization} is complete. 
\qed
\bigskip

 \noindent \textbf{Acknowledgments}.\label{acknowledgments} The authors would like to thank  Gianni Dal~Maso and Luc Tartar for their friendly help, and Lucio Boccardo, Juan Casado-D\'iaz and Luigi Orsina for having introduced them  to singular semilinear problems. They also would like to thank their own institutions (Dipartimento di Scienze di Base e Applicate per l'Ingegneria, Facolt\`a di Ingegneria Civile e Industriale, Sapienza Universit\`a di Roma, Departamento de Matem\'atica Aplicada y Estad\'istica, Universidad Polit\'ecnica de Cartagena, and Laboratoire  Jacques-Louis Lions, Universit\'e Pierre et Marie Curie Paris VI et  CNRS) for providing the support to reciprocal visits which allowed them to perform the present work. The work of Pedro J. Mart\'inez-Aparicio has been partially supported by the grant MTM2015-68210-P of the Spanish Ministerio de Econom\'ia y Competitividad (MINECO-FEDER), the FQM-116 grant of the Junta de Andaluc\'ia and the grant Programa de Apoyo a la Investigaci\'on de la Fundaci\'on S\'eneca-Agencia de Ciencia y Tecnolog\'ia de la Regi\'on de Murcia 19461/PI/14.

\end{document}